\input amstex
\input amsppt.sty
\loadbold

\TagsOnRight \NoBlackBoxes

\define\Y{\Bbb Y}
\define\Z{\Bbb Z}
\define\C{\Bbb C}
\define\R{\Bbb R}

\define\al{\alpha}
\define\be{\beta}
\define\ga{\gamma}
\define\Ga{\Gamma}
\define\de{\delta}
\define\De{\Delta}

\define\La{\Lambda}
\define\la{\lambda}
\define\si{\sigma}
\define\th{\theta}
\define\epsi{\varepsilon}
\define\om{\omega}
\define\Om{\Omega}

\define\wt{\widetilde}
\define\tht{\thetag}

\define\Prob{\operatorname{Prob}}
\define\supp{\operatorname{supp}}

\define\Fun{{\operatorname{Fun}}}

\define\FS{F\!S}

\define\Sym{\La}
\define\bSym{\Sym^\circ}

\define\up{\uparrow}
\define\down{\downarrow}

\define\ome{\boldsymbol\omega}

\define\pd{\partial}

\topmatter

\title Infinite--dimensional diffusions as limits of random walks on partitions
\endtitle
\author Alexei Borodin and Grigori Olshanski \endauthor

\rightheadtext{Diffusion processes}

\abstract Starting with finite Markov chains on partitions of a natural number
$n$ we construct, via a scaling limit transition as $n\to\infty$, a family of
infinite--dimensional diffusion processes. The limit processes are ergodic;
their stationary distributions, the so--called z--measures, appeared earlier in
the problem of harmonic analysis for the infinite symmetric group. The
generators of the processes are explicitly described.

\medskip
\noindent{\smc Keywords:} Diffusion processes; Thoma's simplex; Infinite
symmetric group; Schur functions; z--measures; Dirichlet forms

\endabstract

 \toc \widestnumber\head{8}

\head {} Introduction \endhead

\head 1. The abstract formalism \endhead

\head 2. A toy example: the Pascal triangle \endhead

\head 3. The z--measures
\endhead

\head 4. The ``up'' and ``down'' operators for the z--measures\endhead

\head 5. Convergence of Markov semigroups \endhead

\head 6. The pre--generator $A_{z,z'}$ as a differential operator
\endhead

\head 7. The limit process
\endhead

\head {} References \endhead
\endtoc

\endtopmatter

\document

\head Introduction \endhead

The present paper originated from our previous study of the problem of harmonic
analysis on the infinite symmetric group $S_\infty$. This problem leads to a
family $\{P_z\}$ of probability measures, the z--measures, which depend on the
complex parameter $z\in\C$. The z--measures live on the Thoma simplex, an
infinite--dimensional compact space $\Om$ which is a kind of dual object to the
group $S_\infty$. The aim of the paper is to introduce stochastic dynamics
related to the z--measures. Namely, we construct a family $\{\ome_z(t)\}$ of
diffusion processes in $\Om$ indexed by the same parameter $z\in\C$. The
processes $\{\ome_z(t)\}$ are obtained from certain Markov chains on partitions
of $n$ in a scaling limit as $n\to\infty$. These Markov chains arise in a
natural way, due to the approximation of the group $S_\infty$ by the increasing
chain of the finite symmetric groups $S_n$. Each z--measure $P_z$ serves as a
unique invariant distribution for the corresponding process $\ome_z(t)$, and
the process $\ome_z(t)$ is ergodic with respect to $P_z$. Moreover, $P_z$ is a
symmetrizing measure, so that $\ome_z(t)$ is reversible. We describe the
spectrum of the generator of the process $\ome_z(t)$ and compute the associated
(pre)Dirichlet form.

As was shown in our previous papers, the z--measures and associated random
point processes (that is, measures on infinite point configurations) are
similar to random point processes coming from random matrix ensembles. A
remarkable fact is that this similarity extends to related dynamical (that is,
time--dependent) models. In particular, the dynamical model that we study in
the present paper is similar to time--dependent random point processes (like
the Dyson Brownian motion) coming from random matrices.

Now let us describe our results and their origin in more detail.

\subhead The Thoma simplex\endsubhead This is the subspace
$\Om\subset[0,1]^\infty\times[0,1]^\infty$ consisting of couples $(\al,\be)$,
where $\al$ and $\be$ are two infinite sequences of nonincreasing nonnegative
real numbers, $\al_1\ge\al_2\ge\dots\ge0$ and $\be_1\ge\be_2\ge\dots\ge0$, such
that $\sum\al_i+\sum\be_i\le1$. Define the {\it infinite symmetric group\/}
$S_\infty$ as the union of the chain $S_1\subset S_2\subset S_3\subset\dots$ of
finite symmetric groups. By Thoma's theorem \cite{T}, 1964, the points
$\om=(\al,\be)\in\Om$ parameterize the indecomposable normalized characters of
the group $S_\infty$ (in a hidden form, this remarkable result is contained in
an earlier paper by Edrei \cite{Ed}, 1952). In this sense $\Om$ may be regarded
as a dual object to $S_\infty$.

\subhead Origin of the z--measures\endsubhead The problem of harmonic analysis
for the group $S_\infty$, as stated by Kerov, Olshanski, and Vershik
\cite{KOV1}, consists in decomposing the so--called generalized regular
representations ${\bold T}_z$. These representations depend on an arbitrary
complex number $z$ and form a deformation of the biregular representation of
the group $S_\infty\times S_\infty$ in $\ell^2(S_\infty)$. As shown in
\cite{KOV2}, the decomposition of $\bold T_z$ into irreducibles is uniquely
determined by an equivalence class $\bold P_z$ of Borel measures on $\Om$;
moreover, except the coincidence $\bold P_{z}=\bold P_{\bar z}$, these classes
are mutually singular. \footnote{The latter fact is not too surprising because
the carrying space $\Om$ is infinite dimensional.} The structure of $\bold P_z$
substantially depends on whether $z$ is an integral point or not. The former
case was studied in \cite{KOV2}; it turns out that for each $z\in\Z$, $\bold
P_z$ is supported by a countable union of finite--dimensional faces of the
simplex $\Om$. The latter case, $z\in\C\setminus\Z$, which is of primary
interest for us, is more complex. However, in this case, the construction of
$\bold T_z$ provides a distinguished measure in the class $\bold P_z$, which is
our z--measure $P_z$. \footnote{Actually, the definition of the z--measures can
be extended to a larger set of values of the parameter, see \S3 below.
According to this, starting from \S3 we use, instead of the single subscript
$z$, two subscripts $z$, $z'$.} The z--measures were studied in detail in a
series of our papers, see \cite{BO1--7}.

\subhead Measures on Young diagrams \endsubhead Let $\Y_n$ denote
the set of Young diagrams with $n$ boxes. Recall that there is a
bijective correspondence between the Young diagrams $\la\in\Y_n$ and
the irreducible characters of $S_n$, so that $\Y_n$ may be regarded
as the dual object to $S_n$. The construction of the present paper
relies on the fact that any probability measure $P$ on $\Om$ comes
with a {\it canonical\/} approximation by a sequence $\{M_n\}$ of
probability measures carried by finite sets $\Y_n$. Note that as $n$
gets large, the finite sets $\Y_n$ approximate, in an appropriate
way, the space $\Om$. This intuitively agrees with the fact that the
infinite group $S_\infty$ is approximated by the finite groups
$S_n$. Furthermore, for each $n$, the embedding $S_n\hookrightarrow
S_{n+1}$ induces a canonical Markov transition kernel
$p^{\down}_{n+1,n}$ from $\Y_{n+1}$ to $\Y_n$. It turns out that
$M_n$ and $M_{n+1}$ are always consistent with this kernel: its
application to $M_{n+1}$ gives $M_n$. This can be written as
$M_{n+1}p^\down_{n+1,n}=M_n$; here and below we let a kernel act on
a measure on the right. There exists another transition kernel,
denoted as $p^\up_{n,n+1}$, such that, conversely,
$M_np^\up_{n,n+1}=M_{n+1}$.

This picture can also be described as follows: Any sequence $\{M_n\}$ of
probability measures consistent with the ``down'' transition kernels can be
viewed as the law of a Markov sequence $\{\la^{(n)}\}$ of random variables such
that, for each $n=1,2,\dots$, the $n$th variable $\la^{(n)}$ takes values in
$\Y_n$ and is distributed according to $M_n$, while the ``down'' transition
kernel $p^\down_{n+1,n}$ describes the conditional distribution of $\la^{(n)}$
given $\la^{(n+1)}$. Then the ``up'' transition kernel $p^\up_{n,n+1}$
determines the conditional distribution of $\la^{(n+1)}$ given $\la^{(n)}$.
\footnote{In this context, the Thoma simplex appears as the {\it entrance
boundary\/} for the sequence $\{\la^{(n)}\}$ viewed as a Markov chain with
discrete time $n$ ranging in reverse direction, from $+\infty$ to 1, see
\cite{KOO}. This Markov chain is a little bit unusual and highly
non--homogeneous: its state space varies with time. Nevertheless, the theory of
boundaries (see, e.g., \cite{KSK}) can be readily adapted to such a situation.}

It should be noted that although the ``down'' and ``up'' transition kernels
play symmetric roles, there is a substantial difference between them, for the
``down'' kernel is a canonical object associated with the inductive limit group
$S_\infty=\varinjlim S_n$\,, while the ``up'' kernel varies depending on the
concrete sequence $\{M_n\}$.

\subhead The ``up--down'' Markov chains \endsubhead The starting
point of our construction is the observation that for any $n$ one
can build a ``natural'' Markov chain in $\Y_n$ with stationary
distribution $M_n$. We call it the $n$th {\it up--down chain\/}: By
definition, its transition operator $T_n$  is defined as the
superposition $T_n=p^\up_{n,n+1}\circ p^\down_{n+1,n}$ (the order of
reading is from left to right when $T_n$ is applied to measures, and
from right to left when applied to functions).  Note that if
$(\dots,{\la(t-1)},\la(t),\la(t+1),\dots)$ stands for a trajectory
of our Markov chain in $\Y_n$ then for any time moment $t$, the
Young diagrams $\la(t)$ and $\la(t+1)$ either coincide or differ
from each other by a minimal possible transformation preserving the
total number of boxes: displacement of a single boundary box to a
new position.

One can introduce the graph with the vertex set $\Y_n$ and the edges formed by
couples of Young diagrams such that their symmetric difference consists of two
boxes.  Our chain jumps along the edges and so it is a random walk on this
graph.

Now the idea is to look at the limit behavior of the up--down  Markov chains as
$n\to\infty$. We do not know what can be said in the abstract context, i.e.,
when $\{M_n\}$ comes from an arbitrary probability measure $P$ on $\Om$, but in
the concrete case of the z--measures $P=P_z$ we are able to prove the
convergence of the chains to  continuous time Markov processes in $\Om$ which
turn out to be diffusion processes. Here the limit transition assumes
appropriate scalings both of space (from $\Y_n$ to $\Om$) and of time (from
discrete to continuous). The space scaling is the same as that leading from
$M_n$ to $P$: roughly speaking, it consists in shrinking a diagram $\la\in\Y_n$
by the factor of $n^{-1}$. The time scaling makes one step of the $n$th chain
with large number $n$ equivalent to a small time increment $\Delta t\approx
n^{-2}$.

Note that, instead of the up--down chains, we could equally well
deal with the down--up chains corresponding to the transition
operators of the form $p^\down_{n,n-1}\circ p^\up_{n-1,n}$. At the
finite level one gets slightly different chains but the limit
process does not change.

The above definition of the up--down and down--up Markov chains does
not pretend to originality. The second named author learnt it long
ago from Sergei Kerov (but Kerov never published it). In a different
context, such Markov chains appeared in Fulman's work \cite{F}: he
was interested in the Plancherel distributions $M_n$ which
correspond to the case when $P$ is the delta measure at the
distinguished point $(\al\equiv0,\be\equiv0)$ of $\Om$. Quite
possibly, the trick is well known, for it is very simple and can be
applied to any Markov sequence
$(\dots,\xi^{(n-1)},\xi^{(n)},\xi^{(n+1)},\dots)$ of random
variables with values in possibly varying spaces.

\subhead Limit transition \endsubhead Our analysis of the large $n$ asymptotics
of the up--down Markov chains uses the technique explained in Ethier--Kurtz
\cite{EK2} and goes as follows. Let $C(\Om)$ be the Banach space of continuous
functions with the supremum norm on the compact space $\Om$. We embed $\Y_n$
into $\Om$ by means of a map first introduced by Vershik and Kerov, which gives
us a projection of $C(\Om)$ onto the (finite--dimensional) Banach space
$C(\Y_n)$. The Markov operator $T_n$ can be viewed as an operator in the latter
space. We show that, in an appropriate sense,
$$
\lim_{n\to\infty}n^2(T_n-1)=A, \tag0.1
$$
where $A$ is a closable, dissipative operator $A$ in $C(\Om)$ with dense range.
By the Hille--Yosida theorem and a general theorem due to Kurtz, the closure
$\bar A$ of the operator $A$ generates a contractive semigroup $\{T(t)\}$ in
$C(\Om)$, and moreover, we have convergence of semigroups:
$$
\lim_{n\to\infty}T_n^{[n^2t]}=T(t).
$$
This implies that the limit semigroup $\{T(t)\}$ is positivity preserving.
Consequently, it determines a Feller Markov process in $\Om$, which is our
process $\ome_z(t)$.

The key moment of this argument is the existence of the limit operator $A$ with
desired properties. Surprisingly enough, we can prove this by purely algebraic
tools, without any analytical machinery. Specifically, we are dealing with a
nice algebra of functions on the set $\Y:=\cup\Y_n$ of all Young diagrams. This
algebra, introduced in \cite{KO}, is isomorphic to the algebra $\Sym$ of
symmetric functions, and it turns out that all necessary computations can be
carried out by manipulations with symmetric functions. Here we substantially
use the so--called Frobenius--Schur symmetric functions introduced in
\cite{OlRV} (they are essentially the same as the shifted Schur functions
\cite{OO}). The domain of $A$, $\operatorname{Dom}(A)\subset C(\Om)$, can be
identified with a quotient $\bSym$ of the algebra $\Sym$ by an ideal; this
quotient can be embedded into $C(\Om)$ as a dense subspace. It turns out that
both $A$ and all $T_n$'s can be viewed as operators in $\bSym\subset C(\Om)$
preserving a natural ascending filtration of $\bSym$ by finite--dimensional
subspaces. Due to this fact, the limit transition \tht{0.1} can be understood
at the level of linear algebra: as convergence of operators in these
finite--dimensional subspaces.

\subhead Results about the limit processes \endsubhead We prove that
$\ome_z(t)$ is a diffusion process, that is, a strong Markov process with
continuous sample paths.

We show that $\ome_z(t)$ has the z--measure $P_z$ as a unique stationary
distribution. Moreover $P_z$ is also a symmetrizing measure, so that the
process is reversible.

We describe explicitly the spectrum of the generator, which turns out to be
discrete. Due to existence of a spectral gap, the process is ergodic.

We describe the pre--generator $A=A_z$ in two different ways: (1) as an
operator in the algebra $\Sym^\circ$ acting on the (images of) the Schur
functions and (2) as an infinite--dimensional differential operator in
appropriate coordinates, which we call the moment coordinates. It is worth
noting that these are {\it not\/} the natural coordinates $\al_i$, $\be_i$ in
$\Om$.

We represent the quadratic form associated to $A_z$ as an integral against
$P_z$ of a ``carr\'e du champs'', which does not depend on the parameter and
admits a simple expression in the moment coordinates.

\subhead Lifting \endsubhead Let $\wt\Om$ denote the cone over $\Om$. Each of
the processes $\ome_z(t)$ admits a ``lifting'' to the cone $\wt\Om$. The result
is a Markov process $\widetilde{\ome}_z(t)$ in $\wt\Om$ which is related to
$\ome_z(t)$ by a skew product construction, like the Brownian motions in the
Euclidean space and in its unit sphere. It turns out that, as an equilibrium
process, $\widetilde{\ome}_z(t)$ is a time dependent {\it determinantal
process\/}: its dynamical (i.e., space--time) correlation functions are
principal minors of a kernel. That kernel appeared in our paper \cite{BO6} as
the result of a scaling limit transition from some continuous time jump Markov
processes on Young diagrams. The paper \cite{BO6} left open the question
whether the Markov property persists in that limit transition; now we can
settle this question affirmatively. We aim to discuss this issue in a
subsequent paper.

\subhead Diffusions in Kingman's simplex\endsubhead Diffusion
processes in simplices (both in finite and infinite dimensions) were
extensively discussed in the literature in connection to
mathematical models of population genetics. The paper by  Ethier and
Kurtz \cite{EK1} is of special interest to us. These authors studied
diffusions on the so--called Kingman simplex, which can be
identified with the subspace of points in $\Om$ with all
beta--coordinates equal to 0. On Kingman's simplex, there exists a
remarkable family of probability measures, the Poisson--Dirichlet
distributions; the Ethier--Kurtz diffusions preserve these
distributions. As shown in Petrov \cite{P}, the Ethier--Kurtz
diffusions can be constructed by the limit transition from
appropriate up--down Markov chains; moreover, the same approach can
be extended to a larger family of diffusions on Kingman's simplex,
which are associated to Pitman's generalization of the
Poisson--Dirichlet distributions. Note that the original
construction of \cite{EK1} was quite different. As will be shown in
a subsequent paper, both the Ethier--Kurtz diffusions and our
processes $\ome_z(t)$ are particular cases of a more general
construction. This explains the origin of a striking similarity
between the both families. We would like to gratefully acknowledge
the influence of the paper \cite{EK1} (and of the subsequent paper
by Schmuland \cite{S}) on our work.

\subhead Organization of the paper\endsubhead In \S1 we present the
construction of Markov processes from up--down Markov chains in an abstract
form. Our aim here is to single out the formal part of the argument, which can
be applied in other situations, outside the context of the present paper. In
\S2 we illustrate the abstract formalism on a toy example: here the limit
Markov process is a diffusion on $[0,1]$. In \S3 we introduce the necessary
material related to the z--measures. In \S4 we compute the action of the
transition operators $p^\down_{n+1,n}$ and $p^\up_{n,n+1}$ on symmetric
functions realized as functions on Young diagrams. Using this computation, we
show in \S5 that the assumptions of \S1 are verified for the Markov chain
related to the z--measures; this proves the existence of the limit processes
$\ome_z(t)$. In \S6 we derive an alternative expression for the pre--generator,
which is used in the next section. In \S7, we prove that $\ome_z(t)$ has
continuous sample paths and establish other properties of this process.

\subhead Acknowledgement   \endsubhead The present research was supported  by
the CRDF grant RUM1-2622-ST-04 (both authors), by the NSF grants DMS-0402047
and DMS-0707163 (A.~Borodin), and by the RFBR grant 07-01-91209 and SFB 701,
University of Bielefeld (G.~Olshanski). G.~Olshanski is deeply grateful to Yuri
Kondratiev and Michael R\"ockner for hospitality in Bielefeld and fruitful
discussions.

\head 1. The abstract formalism \endhead

Let $L$ be a graded set, that is, $L$ is the disjoint union of subsets $L_n$,
where $n=0,1,2,\dots$.  Elements of $L$ will be denoted by the letters
$\la,\mu,\nu$. If $\la\in L_n$ then we set $|\la|=n$. We assume that $L_0$ is a
singleton and all $L_n$ are finite. \footnote{The latter assumption could be
relaxed but it is sufficient for the purpose of the present work.} We also
assume that we are given a function $p^\down(\la,\mu)$ on $L\times L$ such
that:
\medskip

$\bullet$ $p^\down(\la,\mu)\ge0$ for all $\la,\mu$;

$\bullet$ $p^\down(\la,\mu)$ vanishes unless $|\la|=|\mu|+1$;

$\bullet$ for every fixed $\la$ with $|\la|\ge1$, we have
$$
\sum_{\mu\in L_{|\la|-1}}p^\down(\la,\mu)=1.
$$

\medskip
Thus, the restriction of $p^\down(\la,\mu)$ to $L_n\times L_{n-1}$ is a
stochastic matrix for every $n\ge1$. We may view this matrix as a transition
function from $L_n$ to $L_{n-1}$. For this reason we call $p^\down$ the {\it
``down'' transition function\/}.

A sequence $M=\{M_0,M_1,\dots\}$, where $M_n$ is a probability
measure on $L_n$, will be called a {\it coherent system\/} (of
measures) if for any $n\ge1$, the measures $M_n$ and $M_{n-1}$ are
consistent with the transition function from $L_n$ to $L_{n-1}$:
$$
\sum_{\la\in L_n}M_n(\la)p^\down(\la,\mu)=M_{n-1}(\mu) \qquad\text{for any
$\mu\in L_{n-1}$,} \tag1.1
$$
where $M_n(\la)$ means the measure of the singleton $\{\la\}$.

Let $M=\{M_n\}$ be a coherent system such that $\supp M_n$, the support of
$M_n$, is the whole $L_n$ for each $n$, that is, $M_n(\la)>0$ for all $n$ and
all $\la\in L_n$. Then we define the {\it``up'' transition function\/} $p^\up$
as follows:
$$
p^\up(\la,\nu)=\frac{M_{n+1}(\nu)}{M_n(\la)}\,p^\down(\nu,\la), \qquad n=|\la|,
\quad \nu\in L_{n+1}\,. \tag1.2
$$

Note that $M=\{M_n\}$ determines a family of random variables indexed by
$n=0,1,\dots$ such that the $n$th variable $\la^{(n)}$ takes the values in
$L_n$ and has the law $M_n$, and, moreover,
$$
\Prob\{\la^{(n-1)}=\mu\mid \la^{(n)}=\la\}=p^\down(\la,\mu)
$$
for each $n\ge1$. In these terms, the ``up'' transition function is interpreted
as the conditional probability
$$
p^\up(\la,\nu)=\Prob\{\la^{(n+1)}=\nu\mid \la^{(n)}=\la\}.
$$
This implies that $p^\up$ is a stochastic matrix of format $L_n\times L_{n+1}$
for any $n=0,1,\dots$.

Note also that each couple $(M_n, M_{n+1})$ is consistent with the ``up''
transition function:
$$
\sum_{\la\in L_n}M_n(\la)p^\up(\la,\nu)=M_{n+1}(\nu) \qquad\text{for any
$\nu\in L_{n+1}$.} \tag1.3
$$

\example{Definition 1.1} Let $M=\{M_n\}$ be a coherent system with
$\supp M_n=L_n$. For any $n$ we define a Markov operator $T_n$ on
the set $L_n$ as the composition of the  ``up'' and ``down''
transition functions, from $L_n$ to $L_{n+1}$ and then back to
$L_n$. The matrix of $T_n$ is given by:
$$
T_n(\la,\wt\la)=\sum_{\nu\in L_{n+1}}p^\up(\la,\nu) p^\down(\nu,\wt\la), \qquad
\la,\wt\la\in L_n.
$$
The operator $T_n$ defines a Markov chain on the set $L_n$.  We call this
Markov chain the ($n$th level) {\it up--down chain\/}. \footnote{Likewise,
interchanging the transition functions, one could introduce the down--up
chains. In the concrete situation studied in the present paper, the down--up
chains slightly differ from the up--down ones but have the same limit.}
\endexample

\proclaim{Proposition 1.2} The measure $M_n$ is an invariant distribution for
the $n$th level up--down Markov chain. Moreover, the chain is reversible with
respect to $M_n$.
\endproclaim

\demo{Proof} The first claim is evident from \tht{1.1} and \tht{1.3}. Indeed,
$p^\up_{n,n+1}$ sends $M_n$ to $M_{n+1}$ and then $p^\down_{n+1,n}$ returns
$M_{n+1}$ back to $M_n$. To prove the second claim we have to check that the
$L_n\times L_n$ matrix $M_n(\la)T_n(\la,\wt\la)$ is symmetric. {}From the
definition of $T_n$ and using \tht{1.2} we get
$$
M_n(\la)T_n(\la,\wt\la) =\sum_{\nu\in
L_{n+1}}M_{n+1}(\nu)p^\down_{n+1,n}(\nu,\la)p^\down_{n+1,n}(\nu,\wt\la),
$$
which is symmetric. \qed
\enddemo

We aim to study the behavior of the up--down chains as $n\to\infty$. To do this
we choose a topological space $\bar L$ and embeddings
$\iota_n:L_n\hookrightarrow\bar L$, and we make an appropriate limit transition
inside $\bar L$. As a result we obtain a continuous time Markov process in
$\bar L$.

An excellent reference about limit transitions from Markov chains to continuous
time Markov processes is Ethier--Kurtz's book \cite{EK2}. We use some general
facts from \cite{EK2} but we also explore some specific properties of our
model. \footnote{The motivation for the concrete choice of the data $(L=(L_n),
p^\down, M=\{M_n\}, \bar L, \iota_n)$ comes from Vershik--Kerov's theory
\cite{VK}, \cite{K}.} Namely, in our concrete situation the following
assumptions hold true:

\medskip

(A1) The ambient space $\bar L$ is a compact, metrizable, separable topological
space. Below we denote by $C(\bar L)$ the Banach space of continuous real
functions with the canonical norm
$$
\Vert f\Vert=\sup_{\om\in\bar L}|f(\om)|.
$$

(A2) The sets $\iota_n(L_n)$ approximate $\bar L$ in the following sense: any
open subset of $\bar L$ has a nonempty intersection with $\iota_n(L_n)$ for all
$n$ large enough.

(A3) There is a distinguished dense subspace $\Cal F$ of the Banach space
$C(\bar L)$ and an ascending exhaustive filtration $(\Cal F^m)$ of $\Cal F$ by
finite--dimensional subspaces such that each $\Cal F^m$ is invariant under the
Markov operators $T_n$ in the following sense:

Denote by $C(L_n)$ the (finite--dimensional) Banach space of functions on $L_n$
with the norm
$$
\Vert g\Vert_n=\sup_{\la\in L_n}|g(\la)|
$$
and let $\pi_n:C(\bar L)\to C(L_n)$ be defined by
$$
(\pi_n(f))(\la)=f(\iota_n(\la)), \qquad \la\in L_n\,, \quad f\in C(\bar L).
$$
Observe that $\pi_n$ is injective on $\Cal F^m$ provided that $n$ is large
enough, where the necessary lower bound on $n$ depends on $m$. \footnote{This
is true because of (A2) and the fact that $\dim\Cal F^m<\infty$. } Then the
invariance property means that, for any $m$, $\pi_n(\Cal F^m)$ is invariant
under $T_n$, at least for large enough $n$.

Thus, identifying $\Cal F^m$ and $\pi_n(\Cal F^m)$, which makes sense for large
$n$, one may say that the operators $T_n$ leave the spaces $\Cal F^m$
invariant.

(A4) There is a sequence $\{\epsi_n\}$ of positive numbers converging to 0 such
that, under the identification $\Cal F^m=\pi_n(\Cal F^m)$, the limit
$$
\lim_{n\to\infty}\epsi_n^{-1}(T_n-\bold1)f=Af.
$$
exists in any finite--dimensional space $\Cal F^m$. Here $\bold1$ stands for
the identity operator. Clearly, we get in this way a limit operator $A:\Cal
F\to\Cal F$ which preserves each subspace $\Cal F^m$.

(A5) $\Cal F$ contains the constant function 1.

\medskip

Recall that a {\it conservative Markov semigroup\/} in $C(\bar L)$ is a
strongly continuous semigroup $\{T(t)\}_{t\ge0}$ of contractive operators in
$C(\bar L)$ preserving the cone of nonnegative functions and the constant
function 1.

\example{Definition 1.3} Following \cite{EK2, Ch. 1, Section 6}, let us say
that a sequence of functions $\{f_n\in C(L_n)\}$ converges to a function $f\in
C(\bar L)$ if $\Vert f_n-\pi_n(f)\Vert_n\to0$. Then we write $f_n\to f$.
\endexample

Observe that $\Vert\pi_n(f)\Vert_n\le\Vert f\Vert$ and, by virtue of (A2),
$\Vert\pi_n(f)\Vert_n\to\Vert f\Vert$. Again by (A2), a sequence $\{f_n\in
C(L_n)\}$ may have at most one limit in $C(\bar L)$.

\proclaim{Proposition 1.4} Assume that the assumptions \tht{A1}--\tht{A5}
stated above are satisfied.

\tht{1} The operator $A:\Cal F\to\Cal F$ defined in \tht{A4} is closable in the
Banach space $C(\bar L)$.

\tht{2} The closure $\bar A$ of $A$ generates a conservative Markov semigroup
$\{T(t)\}$ in $C(\bar L)$.

\tht{3}  The discrete semigroup $\{1,T_n,T_n^2,T_n^3, \dots\}$ converges, as
$n\to\infty$, to $\{T(t)\}$ in the following sense:
$$
T_n^{[\epsi_n^{-1}t]}\pi_n(f)\to T(t)f, \qquad \forall f\in C(\bar L), \tag1.4
$$
for all $t\ge0$, uniformly on bounded intervals, where the limit is understood
according to Definition 1.3.
\endproclaim

We will call $A$ the {\it pre--generator\/} of the semigroup $\{T(t)\}$.

\demo{Proof} Step 1. The operator $A:\Cal F\to\Cal F$ is dissipative, that is,
$\Vert (s\bold1-A)f\Vert\ge s\Vert f\Vert$ for any $s\ge0$ and $f\in\Cal F$.

Indeed, fix $m$ so large that $f\in\Cal F^m$. Assuming $n$ large enough, we may
identify $\Cal F^m$ with its image under $\pi_n$ and view $T_n$ as an operator
in $\Cal F_m$. Set $A_n=\epsi_n^{-1}(T_n-\bold1)$. Since $T_n$ is a contraction
with respect to the norm $\Vert\,\cdot\,\Vert_n$, the operator $T_n-\bold1$ is
dissipative with respect to this norm, whence the same holds for $A_n$. Since
$A_n\to A$ in the finite--dimensional space $\Cal F^m$ and since $\Vert
g\Vert_n\to\Vert g\Vert$ for any $g\in\Cal F^m$ we conclude that $A$ is
dissipative.

Step 2. By virtue of step 1, for any $s>0$ and any $m$, the operator
$s\bold1-A$ maps the finite--dimensional subspace $\Cal F^m$ onto itself. Thus,
$(s\bold1-A)\Cal F=\Cal F$.

Step 3. Since $A$ is dissipative (step 1) and its domain $\Cal F$ is dense, $A$
is closable in $C(\bar L)$; moreover, for any $s>0$, the closure of the range
of $s\bold1-A$ coincides with the range of $s\bold1-\bar A$, where $\bar A$
denotes the closure of $A$ (\cite{EK2, Ch. 1, Lemma 2.11}). Applying step 2 we
see that the range of $s\bold1-\bar A$ coincides with the whole space $C(\bar
L)$. Thus, the operator $\bar A$ satisfies the assumptions of the Hille--Yosida
theorem and therefore it generates a strongly continuous contractive semigroup
$\{T(t)\}$ in $C(\bar L)$, see \cite{EK2, Ch. 1, Thm. 2.6}.

Step 4.  It follows from (A5) that the constant function 1 is in the
domain of $\bar{A}$ and $\bar{A}1=0$, whence $T(t)1=1$.

Step 5. All the assumptions of \cite{EK2, Ch. 1, Thm. 6.5} are satisfied.
Namely:

$\bullet$ $\{T(t)\}$ is a strongly continuous semigroup of contractions with
generator $\bar A$;

$\bullet$ the subspace $\Cal F$ is an essential domain for $\bar A$, that is,
the operator $\bar A$ coincides with the closure of its restriction to $\Cal
F$;

$\bullet$ each $T_n$ is a contraction, and we have
$\epsi_n^{-1}(T_n-\bold1)\pi_n(f)\to \bar Af$ for any $f\in\Cal F$ in the sense
of Definition 1.2.

Applying this theorem from \cite{EK2}, we obtain \tht{1.4}.

Step 6.  Finally, \tht{1.4} implies that the operators $T(t)$ preserve
nonnegative functions, because the $T_n$'s possess this property (here we again
use (A2)). Thus, $\{T(t)\}$ is a Markov semigroup. \qed
\enddemo

\proclaim{Proposition 1.5} The semigroup $\{T(t)\}$ constructed in Proposition
1.4 gives rise to a strong Markov process $X(t)$ in $\bar L$. The process has
c\`adl\`ag sample paths and can start from any point or any probability
distribution.
\endproclaim

\demo{Proof}  This is a well--known general fact, see e.g. \cite{EK2, Ch. 4,
Thm. 2.7}.\qed
\enddemo

\proclaim{Proposition 1.6} Assume additionally that the measures $\iota_n(M_n)$
weakly converge to a measure $P$ on $\bar L$. Then $P$ is an invariant
distribution for the process $X(t)$.
\endproclaim

\demo{Proof} This directly follows from \tht{1.4}. Indeed, let
$\langle\,\cdot\,\rangle_P$ or $\langle\,\cdot\,\rangle_{M_n}$ means
expectation with respect to $P$ or $M_n$. The invariance property of $P$ means
that
$$
\langle T(t)f\rangle_P=\langle f\rangle_P, \qquad \forall f\in C(\bar L), \quad
\forall t\ge0.
$$
Since $\iota_n(M_n)$ weakly converges to $P$, this is equivalent to
$$
\lim_{n\to\infty}\langle T(t)f)\rangle_{\iota_n(M_n)}=\lim_{n\to\infty}\langle
f\rangle_{\iota_n(M_n)},
$$
which can be rewritten as
$$
\lim_{n\to\infty}\langle \pi_n(T(t)f)\rangle_{M_n}=\lim_{n\to\infty}\langle
\pi_n(f)\rangle_{M_n}.
$$
By virtue of \tht{1.4}, as $n$ gets large, $\pi_n(T(t)f)$ is close in norm to
$T_n^{[\epsi_n^{-1}t]}\pi_n(f)$, whence the last limit relation is equivalent
to
$$
\lim_{n\to\infty}\langle
T_n^{[\epsi_n^{-1}t]}\pi_n(f)\rangle_{M_n}=\lim_{n\to\infty}\langle
\pi_n(f)\rangle_{M_n},
$$
which holds for trivial reasons, because
$$
\langle T_n^{[\epsi_n^{-1}t]}\pi_n(f)\rangle_{M_n}=\langle
\pi_n(f)\rangle_{M_n},
$$
due to invariance of $M_n$ with respect to $T_n$. \qed
\enddemo

\proclaim{Proposition 1.7} Under the hypothesis of Proposition 1.6, the
pre--generator $A:\Cal F\to\Cal F$ is symmetric with respect to the inner
product
$$
(f,g):=\langle f\cdot g\rangle_P\,.
$$
\endproclaim

\demo{Proof} The argument is similar to that used in Proposition 1.6. Let us
show that
$$
\langle Af\cdot g\rangle_P
=\lim_{n\to\infty}\langle\epsi_n^{-1}(T_n-\bold1)\pi_n(f)\cdot\pi_n(g)\rangle_{M_n}\,,
\qquad f,g\in\Cal F.
$$
By virtue of Proposition 1.2, the right--hand side is symmetric with respect to
$f\leftrightarrow g$, hence the above equality implies the desired symmetry of
the left--hand side. We have
$$
\gather \langle Af\cdot g\rangle_P=\lim_{n\to\infty}\langle Af\cdot
g\rangle_{\iota_n(M_n)}=\lim_{n\to\infty}\langle \pi_n(Af\cdot
g)\rangle_{\iota_n(M_n)}\\ =\lim_{n\to\infty}\langle
\pi_n(Af)\cdot\pi_n(g)\rangle_{\iota_n(M_n)}=
\lim_{n\to\infty}\langle\epsi_n^{-1}(T_n-\bold1)\pi_n(f)\cdot\pi_n(g)\rangle_{M_n}\,,
\endgather
$$
where the last step is justified using \tht{A4}. \qed
\enddemo

 \proclaim{Proposition 1.8} Under the hypothesis of Proposition
1.6, consider $X(t)$ as an equilibrium process with respect to its invariant
distribution $P$. Likewise, consider the up--down Markov chains in equilibrium
with respect to the invariant distributions $M_n$.

Then the finite--dimensional distributions for the $n$th chain converge, as
$n\to\infty$, to the corresponding finite--dimensional distributions of the
process $X(t)$. Here we assume a natural scaling of time: one step of the $n$th
Markov chain corresponds to a small time interval of order $\Delta t=\epsi_n$.
\endproclaim

\demo{Proof} The argument is similar to that used in Proposition 1.6. \qed
\enddemo

\head 2. A toy example: the Pascal triangle \endhead

Here we illustrate the above formalism on a simple example: the Pascal
triangle.

In this example, the set $L$ consists of arbitrary couples $\la=(a,b)$ of
nonnegative integers. The grading is defined as $|(a,b)|=a+b$. Thus, $L_0$
consists of the single point $(0,0)$, $L_1$ consists of two points $(0,1)$ and
$(1,0)$, \dots, $L_n$ consists of $n+1$ points $(0,n)$, \dots, $(n,0)$.

The ``down'' transition function is defined as follows
$$
p^\down((a,b), (a-1,b))=\frac a{a+b}\,, \qquad p^\down((a,b), (a,b-1))=\frac
b{a+b}\,,
$$
with all other transitions being of probability 0.

As $M_n$ we take the uniform measure on the set $L_n$. The coherency condition
\tht{1.1} is immediately checked. Thus, $M=\{M_n\}$ is a coherent system. The
nonzero values of the ``up'' transition function are
$$
p^\up((a,b), (a+1,b))=\frac{a+1}{a+b+2}\,, \qquad p^\up((a,b), (a,b+1))=\frac
{b+1}{a+b+2}\,.
$$

The Markov operator $T_n$ of the up--down Markov chain is given by the
following matrix (we list the nonzero entries only and assume $a+b=n$)
$$
\align T_n((a,b), (a+1,b-1))&=\frac{(a+1)b}{(a+b+2)(a+b+1)}\\ T_n((a,b),
(a-1,b+1))&=\frac{a(b+1)}{(a+b+2)(a+b+1)}\\
T_n((a,b), (a,b))&=\frac{(a+1)^2+(b+1)^2}{(a+b+2)(a+b+1)}\,.
\endalign
$$

The ambient compact space is the closed unit interval $[0,1]$. The embeddings
$L_n\hookrightarrow\bar L$ are defined as $(a,b)\mapsto x$ with $x=\frac
a{a+b}\in[0,1]$.

The dense subspace $\Cal F\subset C(\bar L)=C([0,1])$ is the space of
polynomials with the canonical filtration by degree.

The fulfilment of  (A1), (A2), and (A5) is evident, let us verify (A3) and
(A4).

The nontrivial one--step transitions of our Markov chain $(a,b)\to(a+1,b-1)$
and $(a,b)\to(a-1,b+1)$ turn into $x\to x\pm\Delta x$ with $\Delta x=n^{-1}$,
where we assume $n=a+b$. According to the above formulas for $T_n$,
$$
\align \Prob\{x\to x+\Delta x\}&=\frac{(a+1)b}{(a+b+2)(a+b+1)}
\\
\Prob\{x\to x-\Delta x\}&=\frac{a(b+1)}{(a+b+2)(a+b+1)} \,,
\endalign
$$
which can be rewritten as
$$
\align \Prob\{x\to x+\Delta x\}&=\frac{n^2}{(n+1)(n+2)}\,
x(1-x)+\frac{n}{(n+1)(n+2)}\,(1-x)
\\
\Prob\{x\to x-\Delta x\}&=\frac{n^2}{(n+1)(n+2)}\,
x(1-x)+\frac{n}{(n+1)(n+2)}\,x \,.
\endalign
$$
It follows that for a polynomial $f(x)$,
$$
\gather ((T_n-\bold1)f)(x)=\frac{n^2}{(n+1)(n+2)}\, x(1-x)(f(x+\Delta
x)+f(x-\Delta x)-2f(x))\\
+\frac{n}{(n+1)(n+2)}(1-x)(f(x+\Delta x)-f(x))+\frac{n}{(n+1)(n+2)}x(f(x-\Delta
x)-f(x))
\endgather
$$
(recall that $\bold1$ stands for the identity operator). Therefore,
$$
\gather ((T_n-\bold1)f)(x)=\frac{n^2\cdot(\De x)^2}{(n+1)(n+2)}\,
x(1-x)(f''(x)+\dots)\\
+\frac{n\cdot\Delta x}{(n+1)(n+2)}(1-x)(f'(x)+\dots) + \frac{n\cdot\Delta
x}{(n+1)(n+2)}x(-f'(x)+\dots),
\endgather
$$
where dots mean higher derivatives multiplied by suitable nonzero
powers of $\Delta x$.

From this expression we see that the operator $T_n$ is well defined
on polynomials and it does not raise the degree, so that assumption
(A3) holds true.

Next, we see that if $\epsi_n\sim n^{-2}$ then, as $n\to\infty$,
$$
\epsi_n^{-1}((T_n-\bold1)f)(x)\quad \to \quad x(1-x)f''(x)+(1-2x)f'(x).
$$
Therefore, assumption (A4) holds with
$$
A=x(1-x)\frac{d^2}{dx^2}+(1-2x)\frac{d}{dx}\,.
$$

Thus, all necessary assumptions are satisfied and one can apply Proposition 1.4
to conclude that our Markov chains converge to a continuous time Markov process
$X(t)$ on $[0,1]$. The generator $\bar A$ of $X(t)$ is the closure of the
differential operator $A$ initially defined on polynomials.

\head 3. The z--measures \endhead

Here we specify the abstract data described in \S1 and introduce related extra
notation. For more detail, see \cite{Ol1}, \cite{KOV2}, \cite{KOO}, \cite{BO3}.

\subhead 3.1. Young diagrams and modified Frobenius coordinates \endsubhead As
$L$ we take the set $\Y$ of all Young diagrams including the empty diagram
$\varnothing$. The subset $L_n\subset L$ becomes the subset $\Y_n\subset\Y$ of
diagrams with $n$ boxes.

Given $\la\in\Y_n$, denote by $a_1,\dots,a_d,b_1,\dots,b_d$ its {\it modified
Frobenius coordinates\/}: here $d$ is the number of diagonal boxes in $\la$,
$a_i$ equals $\frac12$ plus the number of boxes in the $i$th row to the right
of the diagonal, and $b_i$ equals $\frac12$ plus the number of boxes in the
$i$th column below the diagonal. Note that $\sum (a_i+b_i)=n$. We write
$\la=(a_1,\dots,a_d\mid b_1,\dots,b_d)$.

If $\la$ and $\mu$ are two Young diagrams then we write $\mu\nearrow\la$ or,
equivalently, $\la\searrow\mu$ if $\mu\subset\la$ and $|\la|=|\mu|+1$. That is,
$\la$ is obtained from $\mu$ by adding a box. This box is then denoted as
$\la/\mu$. In terms of the modified Frobenius coordinates, $\mu\nearrow\la$
means that $\la$ is obtained from $\mu$ either by adding 1 to one of the
coordinates or by creating a new pair of coordinates $(\frac12;\frac12)$ (the
latter happens if the new box $\la/\mu$ lies on the diagonal).

More generally, for any $\mu\subset\la$ we denote by $\la/\mu$ the
corresponding skew Young diagram.

\subhead 3.2. The ``down'' transition functions \endsubhead The choice of the
``down'' transition function $p^{\down}$ is motivated by the representation
theory of the symmetric groups $S_n$. Recall that $\Y_n$ is the set of labels
of irreducible representations of $S_n$. Given $\la\in\Y_n$, we denote by
$\pi_\la$ the corresponding irreducible representation of $S_n$ and we write
$\dim\la=\dim\pi_\la$. Here is an explicit expression for this quantity in
terms of the modified Frobenius coordinates:
$$
\frac{\dim\la}{n!}=\frac{\prod\limits_{1\le i<j\le d}(a_i-a_j)(b_i-b_j)}
{\prod\limits_{1\le i,j\le d}(a_i+b_j)\prod\limits_{1\le i\le
d}(a_i-\frac12)!(b_i-\frac12)!}\,.
$$

We realize $S_n$ as the group of permutations of the set $\{1,\dots,n\}$, and
we embed $S_{n-1}$ into $S_n$ as the subgroup fixing the point $n$. The Young
rule says that the restriction of $\pi_\la$ (where $\la\in\Y_n$) to
$S_{n-1}\subset S_n$ decomposes into the multiplicity free direct sum of the
representations $\pi_\mu$ such that $\mu\nearrow\la$. Consequently,
$$
\dim\la=\sum_{\mu:\, \mu\nearrow\la}\dim\mu, \qquad |\la|\ge2.
$$
We use this identity to define $p^{\down}$:
$$
p^\down(\la,\mu)=\cases \dfrac{\dim\mu}{\dim\la}\,, & \mu\nearrow\la\\
0, & \text{otherwise}
\endcases
$$
and we also set $p^\down(\la,\mu)=1$ when $|\la|=1$ and $|\mu|=0$, that is,
when $\la$ consists of a single box and $\mu$ is empty.

\subhead 3.3. The Thoma simplex \endsubhead  As $\bar L$ we take the {\it Thoma
simplex\/}. Recall that this is the subspace
$\Om\subset[0,1]^\infty\times[0,1]^\infty$ formed by couples $\om=(\al,\be)$
such that
$$
\gather \al=(\al_1\ge\al_2\ge\dots\ge0)\in[0,1]^\infty, \quad
\be=(\be_1\ge\be_2\ge\dots\ge0)\in[0,1]^\infty, \\
\sum_{i=1}^\infty
(\al_i+\be_i)\le1.
\endgather
$$

The embedding $\iota_n:\Y_n\hookrightarrow\Om$ is defined as follows. For
$\la=(a_1,\dots,a_d\mid b_1,\dots,b_d)\in\Y_n$ (here we wrote $\la$ in terms of
the modified Frobenius coordinates), its image $\iota_n(\la)=(\al,\be)$ is
given by
$$
\al_i=\cases a_i/n, & 1\le i\le d,\\
0, &i>d; \endcases \qquad
\be_i=\cases b_i/n, & 1\le i\le d,\\
0, &i>d.\endcases
$$

The embeddings $\iota_n$ satisfy the assumption (A2) of \S1.

\subhead 3.4. Thoma's measures and moment coordinates \endsubhead To any point
$\om=(\al,\be)\in\Om$ one can assign a probability measure $\nu_\om$ on the
closed interval $[-1,1]$:
$$
\nu_\om=\sum_{i=1}^\infty\al_i\de_{\al_i}+ \sum_{i=1}^\infty\be_i\de_{-\be_i} +
\ga\de_0, \qquad \ga:=1-\sum\al_i-\sum\be_i\,,
$$
where $\de_x$ denotes the Dirac measure at $x$. The measure $\nu_\om$ is called
the {\it Thoma measure\/} corresponding to $\om$.

Denote by $q_k=q_k(\om)$ the moments of $\nu_\om$:
$$
q_k:=\int x^k \nu_\om(dx)=\sum_{i=1}^\infty
\al_i^{k+1}+(-1)^k\sum_{i=1}^\infty\be_i^{k+1}, \qquad k=1,2,\dots
$$
and note that the 0th moment is always equal to 1. We call $q_1, q_2,\dots$ the
{\it moment coordinates\/} of $\om$. Observe that they are continuous functions
in $\om$. Indeed, since $\alpha_i$'s decrease, the condition $\sum\alpha_i\le
1$ implies $\alpha_i\le i^{-1}$ for any $i=1,2,\dots$, whence
$\alpha_i^{k+1}\le i^{-k-1}$. Similarly, $\beta_i^{k+1}\le i^{-k-1}$. It
follows that the both series are uniformly convergent in $\om\in\Om$, which
implies their continuity as functions on $\Om$. \footnote{This argument
substantially relies on the fact that $k+1\ge2$. Note that the function
$\om\mapsto \sum\al_i+\sum\be_i$ is not continuous on $\Om$.}

Let $\Cal M_1[-1,1]$ denote the space of probability Borel measures on $[-1,1]$
equipped with the weak topology. Since this topology is determined by
convergence of moments, the assignment $\om\mapsto\nu_\om$ determines a
homeomorphism of the Thoma simplex on a compact subset of $\Cal M_1[-1,1]$.

On the other hand, the assignment $\om\mapsto (q_1,q_2,\dots)$ determines a
homeomorphism of $\Om$ on a compact subset of $[-1,1]^\infty$.

Thinking of $\Om$ as of a subspace of $\Cal M_1[-1,1]$ or $[-1,1]^\infty$ turns
out to be useful even though we cannot describe the image of $\Om$ in
$[-1,1]^\infty$ explicitly.

Note also that the moment coordinates are algebraically independent as
functions on $\Om$. Indeed, this holds even we restrict them on the subset with
all $\be_i$'s equal to 0. It follows that the algebra of polynomials
$\R[q_1,q_2,\dots]$ can be viewed as a subalgebra of $C(\Om)$, the (real)
Banach algebra of continuous functions on $\Om$ with pointwise operations and
the supremum norm. Since this subalgebra separates points, it is dense in
$C(\Om)$.

\subhead 3.5. Symmetric functions
\endsubhead Let $\Sym$ be the {\it algebra of symmetric functions\/} \cite{Ma}. Recall
that $\Sym$ is freely generated (as a commutative unital algebra) by the {\it
Newton power sums\/} $p_1,p_2,\dots$. As the base field, it is convenient for
us to take $\R$. A distinguished basis in $\La$ is formed by the {\it Schur
functions\/} $s_\mu$ (here and below $\mu$ ranges over $\Y$).

Let $\bSym=\Sym/(p_1-1)$ be the quotient of the algebra $\Sym$ by the ideal
generated by $p_1-1$. The algebra $\Sym^\circ$ has a natural structure of a
filtered algebra inherited from $\Sym$, and the graded algebra associated to
$\bSym$ is isomorphic to $\Sym/(p_1)$.

Given $f\in\Sym$ we denote its image in $\bSym$ by $f^\circ$. In particular, we
will be dealing with the elements $s^\circ_\mu\in\bSym$ coming from the Schur
functions and the elements $p^\circ_k\in\bSym$ coming from the Newton power
sums. Clearly, $p^\circ_1=1$ and $\bSym$ is freely generated (as a unital
commutative algebra) by $p^\circ_2,p^\circ_3, \dots$. \footnote{One could
identify $\bSym$ with the {\it subalgebra\/} in $\Sym$ generated by
$p_2,p_3,\dots$ but we do not want to do this.}

Setting $p^\circ_2\to q_1$, $p^\circ_3\to q_2$, \dots, where $q_k=q_k(\om)$ are
the moment coordinates defined above, we define an algebra isomorphism between
$\bSym$ and the subalgebra $\R[q_1,q_2,\dots]\subset C(\Om)$. Thus, each
element $f^\circ\in\bSym$ becomes a continuous function $f^\circ(\om)$ on
$\Om$, in particular,
$$
p^\circ_k(\om)=\sum_{i=1}^\infty \al_i^k+(-1)^{k-1}\sum_{i=1}^\infty\be_i^k,
\qquad \om=(\al,\be)\in\Om, \quad k=2,3,\dots.
$$
Note that $p^\circ_1\equiv1$.

We will take $\bSym$ as the dense subspace in $C(\Om)$ required in the
assumption (A3) of \S1.

\subhead 3.6. Boundary measures \endsubhead  Let $\{M_n\}$ be an arbitrary
coherent system of probability measures on the sets $\Y_n$ with respect to the
``down'' transition functions introduced in \S3.2;  see \S1 for the general
definition of coherent systems. Denote by $\iota_n(M_n)$ the measure on $\Om$
obtained as the push--forward of $M_n$ with respect to the embedding
$\iota_n:\Y_n\to\Om$ defined in \S3.3.

\proclaim{Theorem} There exists a weak limit $P=\lim_{n\to\infty}\iota_n(M_n)$
on $\Om$. Conversely, $\{M_n\}$ can be reconstructed from $P$ by means of the
equation
$$
M_n(\nu)=\dim\nu\,\int_\Om s^\circ_\nu(\om)P(d\om), \qquad \nu\in\Y_n\,.
$$
The correspondence $\{M_n\}\to P$ is a bijection between coherent systems on
$\Y=\cup\Y_n$ and probability measures on $\Om$.
\endproclaim

Recall that $s^\circ_\nu$ stands for the image in $\La^\circ$ of the Schur
function $s_\nu\in\La$. Since $\La^\circ$ is embedded in $C(\Om)$, the value
$s^\circ_\nu(\om)$ at a point $\om\in\Om$ is well defined.

We call $P$ the {\it boundary measure\/} of the coherent system $\{M_n\}$.

This fundamental result is a refinement of Thoma's theorem \cite{T}.
It is essentially due to Vershik and Kerov \cite{VK}, \cite{K}. See
\cite{KOO} for a detailed proof.

\subhead 3.7. The z--measures \endsubhead Now we proceed to the definition of a
distinguished family of coherent systems. Introduce the notation
$$
(z)_\la=\prod_{(i,j)\in\la}z+j-i, \qquad \la\in\Y, \quad z\in\C,
$$
where ``$(i,j)\in\la$'' means that the product is taken over the boxes of
$\la$; here and below we denote by $(i,j)$ the box with row number $i$ and
column number $j$. The difference $j-i$ is called the {\it content\/} of a
box $(i,j)$.

We define likewise $(z)_{\la/\mu}$ for skew diagrams $\la/\mu$ (then the
product is taken over the boxes in $\la/\mu$). This is a generalization of
the Pochhammer symbol
$$
(z)_n=z(z+1)\dots(z+n-1)
$$
which is obtained in the particular case when  $\la$ is $(n)$, the one--row
diagram with $n$ boxes.

Let $z$ and $z'$ be complex numbers such that $zz'\notin\{0,-1,-2,\dots\}$. The
{\it z--measure\/} on the finite set $\Y_n$ is the complex measure
$M^{(n)}_{z,z'}$ with the weights
$$
M^{(n)}_{z,z'}(\la)=\dfrac{(z)_\la(z')_\la}{(zz')_n}\,
\frac{(\dim\la)^2}{n!}\,, \qquad \la\in\Y_n\,.
$$

The z--measure does not change under transposition $z\leftrightarrow z'$. Thus,
instead of $z$ and $z'$, one can also take as parameters $zz'$ and $z+z'$.

It is known that the weights sum to 1,
$$
\sum_{\la\in\Y_n}M^{(n)}_{z,z'}(\la)=1,
$$
and that the z--measures satisfy \tht{1.1}:
$$
\sum_{\la\searrow\mu}M^{(n)}_{z,z'}(\la)p^\down(\la,\mu)
=M^{(n-1)}_{z,z'}(\mu).
$$
This follows from the representation--theoretical construction of \cite{KOV2}.
For a direct proof, see, e.g., \cite{Ol1}.

The weights $M^{(n)}_{z,z'}(\la)$ are strictly positive for all $n$ if and only
if the couple $(z,z')$ belongs to one of the following two sets in $\C^2$:
\medskip

$\bullet$ {\it Principal series\/}: Both $z$ and $z'$ are not real and are
conjugate to each other.

$\bullet$ {\it Complementary series\/}: Both $z$ and $z'$ are real and are
contained in the same open interval of the form $(N,N+1)$, where $N\in\Z$.

\medskip

The union of these two sets admits a nice description in terms of the
coordinates
$$
x=\frac{z+z'}2, \qquad y=zz'-x^2=-\left(\frac{z-z'}2\right)^2.
$$
Namely, this is the domain $\Cal D$ in the real $(x,y)$--plane bounded from
below by the piecewise smooth curve $C$ built from countably many smooth arcs
$C_N$: here $N$ ranges over $\Z$ and each $C_N$ is an arc of a parabola:
$y=-(x-N)^2$, $|x-N|\le\frac12$. The principal series and the complementary
series are described by points $(x,y)\in \Cal D$ with $y>0$ and $y\le0$,
respectively.

Thus, to each value $(z,z')$ of the principal/complementary series
(equivalently, to each $(x,y)\in \Cal D$), a coherent family
$\{M_n=M^{(n)}_{z,z'}\}$ of probability measures is attached. Note that the
support of $M_n$ is the whole set $\Y_n$.

The ``up'' transition function for the z--measures looks as follows: for
$\la\in\Y_n$
$$
p^\up_{z,z'}(\la,\la^\bullet)=\cases
\dfrac{(z)_{\la^\bullet/\la}(z')_{\la^\bullet/\la}}{zz'+n}\,
\dfrac{\dim\la^\bullet}{(n+1)\dim\la}, & \la^\bullet\searrow\la, \\
0, & \text{otherwise.}
\endcases
$$
We will often use the symbol $\la^\bullet$ to denote a diagram $\nu$ such that
$\nu\searrow\la$. Likewise $\la_\bullet$ will denote a diagram $\mu$ such that
$\mu\nearrow\la$. Note that $\la^\bullet/\la$ is a single box.

\subhead 3.8. The boundary z--measures \endsubhead Given a coherent system
$\{M^{(n)}_{z,z'}\}$ of z--measures, we denote by $P_{z,z'}$ the corresponding
boundary measure on $\Om$ and call it the {\it boundary z--measure\/}.

\proclaim{Theorem} Except the equality $P_{z,z'}=P_{z',z}$, the boundary
z--measures with different parameters are mutually singular with respect to
each over.
\endproclaim

See \cite{KOV2} for a proof.

\head 4. The ``up'' and ``down'' operators for the z--measures \endhead

Let $\Fun(\Y)$ denote the algebra of all real--valued functions on $\Y$ with
pointwise operations. We define an algebra morphism $\Sym\to\Fun(\Y)$ by
specifying it on the generators $p_k$, as follows
$$
p_k(\la)=\sum_{i=1}^d  a_i^k+(-1)^{k-1}\sum_{i=1}^d b_i^k, \qquad \la\in\Y,
$$
where $a_1,\dots,a_d$ and $b_1,\dots,b_d$ are the modified Frobenius
coordinates of $\la$. It is readily verified that the images of the
$p_k$'s in $\Fun(\Y)$ are algebraically independent, so that our
morphism is injective. Thus, we may view $\Sym$ as a subalgebra of
$\Fun(\Y)$. We will denote the value at $\la$ of the function on
$\Y$ corresponding to an element $f\in\Sym$ as $f(\la)$.

As mentioned above, the Schur functions $s_\mu$ form a distinguished basis in
$\La$. The Schur functions are homogeneous, the degree of $s_\mu$ is equal to
$|\mu|$, see \cite{Ma}.

There is another important basis in $\La$, formed by the {\it
Frobenius--Schur\/} functions $\FS_\mu$. These are inhomogeneous elements such
that $\FS_\mu$ differs from $s_\mu$ by lower degree terms. The crucial property
of the $\FS_\mu$'s is expressed by the formula
$$
\FS_\mu(\la)=n^{\down m}\,\frac{\dim(\mu,\la)}{\dim\la}\,, \qquad n=|\la|,
\quad m=|\mu|, \tag4.1
$$
where we use the notation
$$
n^{\down m}=n(n-1)\dots(n-m+1)
$$
and $\dim(\mu,\la)$ denotes the number of all possible chains
$\mu\nearrow\dots\nearrow\la$ leading from $\mu$ to $\la$ (equivalently,
$\dim(\mu,\la)$ equals the number of standard tableaux of the skew shape
$\la/\mu$ if $\la$ contains $\mu$, and 0 otherwise). In particular,
$\FS_\mu(\la)$ vanishes unless $\la$ contains $\mu$.

For more detail about the realization of symmetric functions as functions on
$\Y$ and about the Frobenius--Schur functions, see \cite{KO}, \cite{OO},
\cite{OlRV}, \cite{IO}.

Given $f\in\Sym$, we denote by $f_n$ the restriction of the function
$f(\,\cdot\,)$ to $\Y_n\subset\Y$. It is readily checked that the subalgebra
$\Sym\subset\Fun(\Y)$ separates points, which implies that for each $n$, the
functions of the form $f_n$, with $f\in\Sym$, exhaust the space $C(\Y_n)$.

Let $D_{n+1,n}: C(\Y_n)\to C(\Y_{n+1})$ and $U_{n,n+1}: C(\Y_{n+1})\to C(\Y_n)$
be the ``down'' and ``up'' operators acting on functions:
$$
\gather (D_{n+1,n}f)(\nu)=\sum_{\la\in\Y_n}p^\down(\nu,\la)f(\la), \qquad f\in
C(\Y_n), \quad
\nu\in\Y_{n+1}\,,\\
(U_{n,n+1}g)(\la)=\sum_{\nu\in\Y_{n+1}}p_{z,z'}^\up(\la,\nu)g(\nu), \qquad g\in
C(\Y_{n+1}), \quad\la\in\Y_n\,.
\endgather
$$
Note that $U_{n,n+1}$ depends on $z$ and $z'$ while $D_{n+1,n}$ does not.

In this section, we prove the following claim.

\proclaim{Theorem 4.1} \tht{1} There exists a unique operator $\wt
D:\Sym\to\Sym$ such that
$$
D_{n+1,n}f_n=\frac1{n+1}(\wt Df)_{n+1}\,, \qquad\text{for all $n=0,1,\dots$ and
all $f\in\Sym$}.
$$
In the basis $\{\FS_\mu\}$ it is given by
$$
\wt D\, \FS_\mu=(p_1-|\mu|)\FS_\mu\,, \qquad \mu\in\Y.
$$

\tht{2} There exists a unique operator $\wt U:\Sym\to\Sym$ depending on $z,z'$,
such that
$$
U_{n,n+1}f_{n+1}=\frac1{zz'+n}(\wt Uf)_n\,, \qquad\text{for all $n=0,1,\dots$
and all $f\in\Sym$}.
$$
In the basis $\{\FS_\mu\}$ it is given by
$$
\wt U\,
\FS_\mu=\sum_{\mu_\bullet\nearrow\mu}(z)_{\mu/\mu_\bullet}(z')_{\mu/\mu_\bullet}
\FS_{\mu_\bullet}+(p_1+zz'+|\mu|)\FS_\mu\,,\qquad \mu\in\Y.
$$
\endproclaim

\demo{Proof of\/ \tht{1}} Uniqueness follows from the fact that
$\Sym\to\Fun(\Y)$ is an embedding. Let us check the required relation. Denote
$m=|\mu|$. We have to prove that for any $n$ and any $\nu\in\Y_{n+1}$,
$$
(D_{n+1,n}(\FS_\mu)_n)(\nu)=\frac{n+1-m}{n+1}\FS_\mu(\nu)
$$
(here we have used the fact that $p_1(\nu)=n+1$).

If $n<m$ then the equality holds for trivial reasons: both sides vanish.
Indeed, we have $(\FS_\mu)_n\equiv0$, the factor $n+1-m$ vanishes for $n=m-1$,
and $\FS_\mu(\nu)=0$ for $n<m-1$. Thus, we may assume $n\ge m$.

Then we have
$$
\align
(D_{n+1,n}(\FS_\mu)_n)(\nu)&=\sum_{\nu_\bullet\nearrow\nu}\frac{\dim\nu_\bullet}
{\dim\nu}\, \FS_\mu(\nu_\bullet)\\
&=\sum_{\nu_\bullet\nearrow\nu}\frac{\dim\nu_\bullet} {\dim\nu}\,n^{\down
m}\,\frac{\dim(\mu,\nu_\bullet)}{\dim\nu_\bullet} \qquad\text{by \tht{4.1}}\\
&=n^{\down m}\sum_{\nu_\bullet\nearrow\nu} \frac{\dim(\mu,\nu_\bullet)}{\dim\nu}\\
&=n^{\down m}\frac{\dim(\mu,\nu)}{\dim\nu}\\
&=\frac{n^{\down m}}{(n+1)^{\down m}}\, \FS_\mu(\nu) \qquad\text{by
\tht{4.1}}\\
&=\frac{n+1-m}{n+1}\, \FS_\mu(\nu),
\endalign
$$
as required. \qed
\enddemo

The proof of \tht{2} is more involved and depends on the lemma below which is
essentially due to Sergei Kerov (see Okounkov \cite{Ok}).

Let $\Fun_0(\Y)\subset\Fun(\Y)$ be the space of functions with finite support,
and let $\{\de_\la\}$ be its natural basis: $\de_\la(\la)=1$ and
$\de_\la(\nu)=0$ for $\nu\ne\la$. Consider the Lie algebra $\frak{sl}(2,\C)$
with its basis
$$
E=\bmatrix 0 & 1\\ 0 & 0 \endbmatrix, \quad F=\bmatrix 0 & 0\\ 1 & 0
\endbmatrix, \quad H=\bmatrix 1 & 0\\ 0 & -1 \endbmatrix.
$$

\proclaim{Lemma 4.2} For any complex $z$ and $z'$, the following action of $E$,
$F$, and $H$ in the basis $\{\de_\la\}$ defines a representation of
$\frak{sl}(2,\C)$ in $\Fun_0(\Y)$
$$
E\de_\la=\sum_{\la^\bullet\searrow\la}
(z)_{\la^\bullet/\la}(z')_{\la^\bullet/\la}\de_{\la^\bullet}\,, \quad
F\de_\la=-\sum_{\la_\bullet\nearrow\la}\de_{\la_\bullet}\,, \quad
H\de_\la=(zz'+2|\la|)\de_\la\,.
$$
\endproclaim
\demo{Proof} The only nontrivial commutation relation to be checked is
$[E,F]=H$. We have
$$
[E,F]\de_\la
=\sum_{\varkappa\nearrow\la^\bullet}\sum_{\la^\bullet\searrow\la}
(z)_{\la^\bullet/\la}(z')_{\la^\bullet/\la}\de_{\varkappa} -
\sum_{\varkappa\searrow\la_\bullet}\sum_{\la_\bullet\nearrow\la}
(z)_{\varkappa/\la_\bullet}(z')_{\varkappa/\la_\bullet}\de_{\varkappa}
$$
The right--hand side is a linear combination of the vectors $\de_\varkappa$
such that either $\varkappa=\la$ or $\varkappa$ is obtained from $\la$ by
adding a box $\square_1$ and removing another box $\square_2\ne\square_1$. In
the latter case, the coefficient of $\de_\varkappa$ in each double sum equals
$(z)_{\square_1}(z')_{\square_1}$, so that the total coefficient is 0.

Examine now the coefficient of $\de_\la$, which is equal to
$$
\sum_{\la^\bullet\searrow\la}(z)_{\la^\bullet/\la}(z')_{\la^\bullet/\la} -
\sum_{\la_\bullet\nearrow\la}(z)_{\la/\la_\bullet}(z')_{\la/\la_\bullet}
$$
Denoting by $\{x_i\}$ and $\{y_j\}$ the contents of the boxes\,
\footnote{Recall that the content of a box $(i,j)$ is defined as $j-i$.} that
can be added to the diagram $\la$ or removed from it, respectively, we write
the above expression as
$$
\multline \sum_i(z+x_i)(z'+x_i) -\sum_j(z+y_j)(z'+y_j)\\
=\left(\sum_i 1-\sum_j 1\right)zz'+\left(\sum_i x_i-\sum_j
y_j\right)(z+z')+\left(\sum_i x_i^2-\sum_j y_j^2\right).
\endmultline
$$

As was first observed by Kerov (see his book \cite{K, Ch. IV, \S1}), the
$x_i$'s and the $y_j$'s form two interlacing sequences
$$
x_1<y_1<x_2<\dots<x_k<y_k<x_{k+1}
$$
such that
$$
\sum_i x_i-\sum_j y_j=0, \qquad \sum_i x_i^2-\sum_j y_j^2=2|\la|.
$$
The easiest way to prove this is to proceed by induction on $|\la|$, by
consecutively adding a box to the diagram.

It follows that in our expression, the coefficient of $zz'$ equals 1 (because
the number of $x$'s is greater than the number of $y$'s by 1), that of $z+z'$
equals 0, and the last term equals $2|\la|$. This completes the proof. \qed
\enddemo

The next lemma is a simple observation:

\proclaim{Lemma 4.3} Let\/ $\square_{k,l}$ denote the rectangular diagram with
$k$ rows and $l$ columns, and $V_{k,l}\subset\Fun_0(\Y)$ stand for the finite
dimensional subspace spanned by the basis vectors $\de_\la$ such that
$\la\subseteq\square_{k,l}$.

If $z=k$ and $z'=-l$ then $V_{k,l}$ is $\frak{sl}(2,\C)$--invariant and the
action of $\frak{sl}(2,\C)$ in $V_{k,l}$ lifts to a representation of the group
$SL(2,\C)$.
\endproclaim

\demo{Proof} If a diagram $\la$ is contained in $\square_{k,l}$
while a diagram $\la^\bullet$, such that $\la^\bullet\searrow\la$,
is not, then the square $\la^\bullet/\la$ may be only one of the
boxes $(1,l+1)$ or $(k+1,1)$. In the former case,
$(z')_{\la^\bullet/\la}$ vanishes, and in the latter case
$(z)_{\la^\bullet/\la}$ vanishes. Therefore, the coefficient of
$\de_{\la^\bullet}$ in the expansion of $E\de_\la$ equals 0. It
follows that the subspace $V_{k,l}$ is $E$--invariant, and its
invariance with respect to $F$ and $H$ is obvious. Thus, $V_{k,l}$
is an $\frak{sl}(2,\C)$--module. Since it has finite dimension, it
generates a representation of the group $SL(2,\C)$. Note that this
representation is irreducible but we do not need this fact. \qed
\enddemo

We proceed to the proof of the second claim of Theorem 4.1.

\demo{Proof of\/ \tht{2}} As in (1), the uniqueness part of the claim is
evident. The remaining (nontrivial) part of the claim means that for any $n$
and any $\la\in\Y_n$,
$$
\aligned (zz'+n)(U_{n,n+1}(\FS_\mu)_{n+1})(\la)
&=\sum_{\mu_\bullet\nearrow\mu}(z)_{\mu/\mu_\bullet}(z')_{\mu/\mu_\bullet}
\FS_{\mu_\bullet}(\la)\\
&+(n+zz'+|\mu|)\FS_\mu(\la).
\endaligned \tag4.2
$$
For $\mu=\varnothing$,  $\FS_\mu$ reduces to the constant function 1, the sum
in the right--hand side disappears, and \tht{4.2} reduces to the tautology
$zz'+n=zz'+n$.

Assume now $|\mu|=m\ge1$. Using the definition of $U_{n, n+1}$ and the basic
formula \tht{4.1}, one can reduce \tht{4.2} to the following combinatorial
identity
$$
\aligned
\sum_{\la^\bullet\searrow\la}(z)_{\la^\bullet/\la}(z')_{\la^\bullet/\la}\,
n^{\down(m-1)}\dim(\mu,\la^\bullet)
&=\sum_{\mu_\bullet\nearrow\mu}(z)_{\mu/\mu_\bullet}(z')_{\mu/\mu_\bullet}
n^{\down(m-1)}\dim(\mu_\bullet,\la)\\
&+(n+zz'+m)n^{\down m}\dim(\mu,\la),
\endaligned
$$
If $n<m-1$ then both sides vanish. Thus, we may assume $n\ge m-1$, so that
$n^{\down(m-1)}\ne0$. Dividing by $n^{\down(m-1)}$ we reduce the identity to
$$
\aligned
\sum_{\la^\bullet\searrow\la}(z)_{\la^\bullet/\la}(z')_{\la^\bullet/\la}\,
\dim(\mu,\la^\bullet)
&=\sum_{\mu_\bullet\nearrow\mu}(z)_{\mu/\mu_\bullet}(z')_{\mu/\mu_\bullet}
\dim(\mu_\bullet,\la)\\
&+(n+zz'+m)(n-m+1)\dim(\mu,\la),
\endaligned \tag4.3
$$

Observe that the identity is satisfied if $\la$ does not contain $\mu$. Indeed,
in such a case $\dim(\mu,\la)=0$, and the last summand disappears. If the set
difference $\mu\setminus\la$ contains 2 or more boxes then no $\la^\bullet$
contains $\mu$ and no $\mu_\bullet$ is contained in $\la$, so that both sides
vanish. Examine now the case when $\mu\setminus\la$ consists of a single square
$\square$. Then the only nonzero contribution to the left--hand side comes from
the summand with $\la^\bullet=\la\cup\square$, and the only nonzero
contribution to the right--hand side comes from
$\mu_\bullet=\mu\setminus\square$. Since
$\la^\bullet/\la=\mu/\mu_\bullet=\square$, the identity is reduced to
$\dim(\mu,\la^\bullet)=\dim(\mu_\bullet,\la)$, which is obvious, because the
skew diagrams $\la^\bullet/\mu$ and $\la/\mu_\bullet$ coincide.

Thus, we may assume $\mu\subseteq\la$. We will check the identity using Lemma
4.2. Since both sides are polynomials in $z$ and $z'$, we may assume that $z=k$
and $z'=-l$, where $k$ and $l$ are so large that all diagrams $\la^\bullet$ are
contained in  $\square_{k,l}$.

Let us multiply the left-hand side of \tht{4.3} by
$(z)_{\la/\mu}(z')_{\la/\mu}$. Due to our assumptions this quantity
is well defined and is nonzero. We obtain
$$
\gather (z)_{\la/\mu}(z')_{\la/\mu}
\sum_{\la^\bullet\searrow\la}(z)_{\la^\bullet/\la}(z')_{\la^\bullet/\la}\,
\dim(\mu,\la^\bullet)
=\sum_{\la^\bullet\searrow\la}(z)_{\la^\bullet/\mu}(z')_{\la^\bullet/\mu}\,
\dim(\mu,\la^\bullet)\\
=\sum_{\la^\bullet\searrow\la}(E^{n-m+1}\de_\mu,\de_{\la^\bullet})
=(n+1-m)!\sum_{\la^\bullet\searrow\la}(e^E\de_\mu,\de_{\la^\bullet})\\
=(n+1-m)!(e^E\de_\mu,\sum_{\la^\bullet\searrow\la}\de_{\la^\bullet}),
\endgather
$$
where all operators act in the finite--dimensional subspace $V_{k,l}$ described
in Lemma 4.3, and $(\,\cdot\,,\,\cdot\,)$ is the natural inner product
inherited from $\ell^2(\Y)$.

Since the operator $\de_\la\mapsto
\sum_{\la^\bullet\searrow\la}\de_{\la^\bullet}$ is adjoint to $-F$, our
expression can be rewritten simply as
$$
-(n+1-m)!(Fe^E\de_\mu,\de_\la).
$$

A simple computation in $SL(2,\C)$ shows that
$$
-Fe^E=-e^E(e^{-E}F e^E)=e^E(-F+E+H),
$$
and due to the last claim of Lemma 4.3 we may interpret the above identity as a
relation between operators in $V_{k,l}$. Then we obtain
$$
\gather -(Fe^E\de_\mu,\de_\la)=-(e^E
F\de_\mu,\de_\la)+(e^EH\de_\mu,\de_\la)+(e^EE\de_\mu,\de_\la)\\
=\sum_{\mu_\bullet\nearrow\mu}(e^E\de_{\mu_\bullet},\de_\la)
+(zz'+2m)(e^E\de_\mu,\de_\la)+(Ee^E\de_\mu,\de_\la) \tag4.4
\endgather
$$

It remains to check that multiplying \tht{4.4} by $(n-m+1)!$ gives the
right--hand side of \tht{4.3} multiplied by $(z)_{\la/\mu}(z')_{\la/\mu}$. The
expression \tht{4.4} comprises three terms.

The first term gives
$$
\gather (n-m+1)!\sum_{\mu_\bullet\nearrow\mu}(e^E\de_{\mu_\bullet},\de_\la)
=\sum_{\mu_\bullet\nearrow\mu}
(z)_{\la/\mu_\bullet}(z')_{\la/\mu_\bullet}\dim(\mu_\bullet,\la)\\
=(z)_{\la/\mu}(z')_{\la/\mu} \sum_{\mu_\bullet\nearrow\mu}
(z)_{\mu/\mu_\bullet}(z')_{\mu/\mu_\bullet}\dim(\mu_\bullet,\la).\tag4.5
\endgather
$$

Next, the second term gives
$$
(n-m+1)!(zz'+2m)(e^E\de_\mu,\de_\la)=(z)_{\la/\mu}(z')_{\la/\mu}(zz'+2m)(n-m+1)\dim(\mu,\la),
$$
the third term gives
$$
(n-m+1)!(Ee^E\de_\mu,\de_\la)=(z)_{\la/\mu}(z')_{\la/\mu}(n-m+1)(n-m)\dim(\mu,\la),
$$
and their sum equals
$$
(z)_{\la/\mu}(z')_{\la/\mu}(zz'+n+m)(n-m+1)\dim(\mu,\la).\tag4.6
$$

We see that the sum of \tht{4.5} and \tht{4.6} is indeed equal to the
right--hand side of \tht{4.3} multiplied by $(z)_{\la/\mu}(z')_{\la/\mu}$. \qed
\enddemo

\head 5. Convergence of Markov semigroups \endhead

Fix arbitrary parameters $(z,z')$ of principal or complementary series. Let
$\{M_n\}=\{M^{(n)}_{z,z'}\}$ be the corresponding coherent family of
probability measures on $\Y$. According to the general formalism of \S1, we
form, for each $n=1,2,\dots$, the $n$th level up--down Markov chain on the set
$\Y_n$ of Young diagrams with $n$ boxes. Since the move ``up'' consists in
appending a box to a Young diagram, while the move ``down'' consists in
removing a box, any nontrivial change of our up--down chain under one step
reduces to moving one of the boxes of a Young diagram to a new position.

The measure $M_n$ is an invariant measure of the chain. It is readily seen that
all the states are communicating, so that $M_n$ is a unique invariant
probability measure.

As explained in \S3, we consider the embeddings $\iota_n:\Y_n\to\Om$ determined
by the normalized modified Frobenius coordinates. Let $\pi_n$ be the
corresponding linear map $C(\Om)\to C(\Y_n)$.

As the space $\Cal F\subset C(\Om)$ we take the algebra $\bSym$. It is dense in
$C(\Om)$ by virtue of \cite{KOO, Lemma 5.3}. The filtration in $\bSym$ is
inherited from $\Sym$.

The two limit relations in the claim below are understood in the
sense of \S1.

\proclaim{Theorem 5.1} With these data, all the assumptions of Proposition 1.4
are satisfied provided that the scaling of time is determined by the factors
$\epsi_n\sim n^{-2}$. Thus, denoting by $T_n$ the Markov operator of the $n$th
level Markov chain, we have convergence to a conservative Markov semigroup
$\{T(t)\}$ in the Banach space $C(\Om)$, as in \tht{1.4}{\rm:}
$$
\lim_{n\to\infty}T^{[n^2t]}\pi_n(f)=T(t)f
$$
for any fixed $f\in C(\Om)$ and all $t\ge0$, uniformly on bounded
intervals. Furthermore, the generator of the limit semigroup
$\{T(t)\}$ is the closure of the operator $A$ with domain
$\bSym\subset C(\Om)$, defined by
$$
Af=\lim_{n\to\infty}n^2(T_n-\bold1)f, \qquad f\in\bSym.
$$
\endproclaim

\demo{Proof} By virtue of Proposition 1.4, it suffices to check assumptions
(A1)--(A5) stated before the formulation of the proposition.

The fulfilment of assumptions (A1), (A2), and (A5) is obvious.

We proceed to verifying assumption (A3).

Recall that the Markov operator $T_n: C(\Y_n)\to C(\Y_n)$ is defined as the
composition $U_{n,n+1}\circ D_{n+1,n}$ of the up and down operators. As shown
in \S4, the latter operators are implemented by certain operators in $\Sym$.
Here we interpret $\Sym$ as a subalgebra in $\Fun(\Y)$ and consider the
restriction map $\Fun(\Y)\to C(\Y_n)$ turning elements  $f\in\Sym$ to functions
$f_n$ on $\Y_n$.

\proclaim{Lemma 5.2} Let $\mu$ be a Young diagram and $m=|\mu|$. The operator
$T_n-\bold1$ acts on $(\FS_\mu)_n$ as follows
$$
\multline
(T_n-\bold1)(\FS_\mu)_n = -\,\frac{m(m-1+zz')}{(n+1)(zz'+n)}(\FS_\mu)_n\\
+\frac{n+1-m}{(n+1)(zz'+n)}\sum_{\mu_\bullet\nearrow\mu}
(z)_{\mu/\mu_\bullet}(z')_{\mu/\mu_\bullet}(\FS_{\mu_\bullet})_n \,.
\endmultline
$$
\endproclaim

\demo{Proof} This follows directly from the computation of \S4. Indeed, we
have
$$
\gather T_n(\FS_\mu)_n=U_{n,n+1}D_{n+1,n}(\FS_\mu)_n\,,\\
D_{n+1,n}(\FS_\mu)_n=\frac{n+1-m}{n+1}(\FS_\mu)_{n+1}\,,\\
U_{n,n+1}(\FS_\mu)_{n+1}=\frac1{zz'+n}
\left(\sum_{\mu_\bullet\nearrow\mu}(z)_{\mu/\mu_\bullet}
(z')_{\mu/\mu_\bullet}(\FS_{\mu_\bullet})_n+(n+zz'+m)(\FS_\mu)_n\right),
\endgather
$$
which implies the desired expression. Note that in the last equality we used
the fact that $(p_1 \FS_\mu)_n=n(\FS_\mu)_n$. \qed
\enddemo

\proclaim{Corollary 5.3} $T_n$ preserves the filtration in $\Cal F$.
\endproclaim

Let us verify assumption (A4) with $\epsi_n=n^{-2}$:

Consider the map $\Sym\to C(\Y_n)$, defined as restriction to
$\iota_n(\Y_n)\subset\Om$. We denote it as $ f\mapsto f_{[n]}$ (it should not
be confused with the map $f\mapsto f_n$ introduced in the beginning of \S4!).
By the very definition, if $\la=(a_1,\dots,a_d\mid b_1,\dots,b_d)\in\Y_n$ (the
modified Frobenius coordinates) then
$$
f_{[n]}(\la)=f\left(\frac{a_1}n,\dots,\frac{a_d}n,0,0,\dots;
\frac{b_1}n,\dots,\frac{b_d}n,0,0,\dots\right), \qquad f\in\Sym.
$$
In the notation of \S1,
$$
f_{[n]}=\pi_n(f^\circ), \qquad f\in\Sym.
$$

Let $G:\Sym\to\Sym$ denote the operator acting on the $m$th homogeneous
component of $\Sym$ as multiplication by $m$ ($m=0,1,2,\dots$). According to
this definition, we denote by $s^G$, where $s\ne0$, the automorphism of the
algebra $\Sym$ that reduces to multiplication by $s^m$ on the homogeneous
component of degree $m$. Then we have
$$
f_n(\la)=(n^Gf)_{[n]}(\la), \qquad \la\in\Y_n.
$$
Indeed, it suffices to check this formula for $f=p_k$, and then it follows from
the very definition of the embedding $\La\to\Fun(\Y)$, see the beginning of
\S4.

\proclaim{Lemma 5.4} There exists a linear operator $A_{z,z'}$ in $\Cal
F=\bSym$ which is the limit of the operators $n^2(T_n-\bold1)$ as $n\to\infty$.
Specifically,
$$
A_{z,z'}s^\circ_\mu=\,-\, m(m-1+zz')s^\circ_\mu +\sum_{\mu_\bullet\nearrow\mu}
(z)_{\mu/\mu_\bullet}(z')_{\mu/\mu_\bullet}s^\circ_{\mu_\bullet}\,, \qquad
\mu\in\Y, \quad m:=|\mu|. \tag5.1
$$
\endproclaim

\demo{Proof} Rewrite the claim of Lemma 5.2 where we substitute
$(\FS_\mu)_n=(n^G \FS_\mu)_{[n]}$ and $(\FS_{\mu_\bullet})_n=(n^G
\FS_{\mu_\bullet})_{[n]}$\,:
$$
\multline
(T_n-\bold1)(n^G \FS_\mu)_{[n]}=-\,\frac{m(m-1+zz')}{(n+1)(zz'+n)}(n^G \FS_\mu)_{[n]}\\
+\frac{n+1-m}{(n+1)(zz'+n)}\sum_{\mu_\bullet\nearrow\mu}
(z)_{\mu/\mu_\bullet}(z')_{\mu/\mu_\bullet}(n^G \FS_{\mu_\bullet})_{[n]} \,.
\endmultline
$$

Multiply both sides by $n^2\cdot n^{-m}$ and observe that
$$
\lim_{n\to\infty}n^{-m}n^G \FS_\mu =s_\mu, \qquad \lim_{n\to\infty}n^{-m+1}n^G
\FS_{\mu_\bullet}=s_{{\mu_\bullet}},
$$
because $\FS_\mu$ and $s_\mu$, as well as $\FS_{\mu_\bullet}$ and
$s_{\mu_\bullet}$, differ in lower order terms only. This implies the claim of
the lemma.\qed
\enddemo

Note that for $\mu=\varnothing$, the sum in \tht{5.1} disappears and $m$
vanishes, so that $A_{z,z'}$ sends $s^\circ_{\varnothing}=1$ to 0, as it should
be.

This concludes the proof of Theorem 5.1. \qed
\enddemo

Looking at formula \tht{5.1}, it is not obvious that it defines an
operator in $\bSym$, because the elements $s^\circ_\mu$ are not
linearly independent. Of course, correctness of \tht{5.1} follows
from the computation in the proof of Lemma 5.4. On the other hand,
this also can be proved directly:

\proclaim{Proposition 5.5} The operator  $B_{z,z'}:\Sym\to\Sym$ determined in
the basis of Schur functions by the expression
$$
B_{z,z'}s_\mu=-|\mu|(|\mu|-1+zz')s_\mu+ p_1\sum_{\mu_\bullet\nearrow\mu}
(z)_{\mu/\mu_\bullet}(z')_{\mu/\mu_\bullet}s_{\mu_\bullet} \tag5.2
$$
preserves the principal ideal generated by $p_1-1$, and the reduction of
$B_{z,z'}$ modulo this ideal coincides with the expression \tht{5.1}.
\endproclaim

\demo{Proof} The second claim is obvious from the comparison of \tht{5.2} with
\tht{5.1}; note that the prefactor $p_1$ in front of the sum in \tht{5.2} will
disappear after the reduction. Note also that the operator $B$ preserves the
grading in $\Sym$.

To prove the first claim it suffices to check that $B$ commutes with the
operator of multiplication by $p_1$.

We use the representation of the Lie algebra $\frak{sl}(2,\C)$ in the space
$\Fun_0(\Y)$ defined in Lemma 4.2. Let $H^*$, $E^*$, and $F^*$ be the adjoint
operators to $H$, $E$, and $F$, respectively. We interpret them again as
operators in $\Fun_0(\Y)$, given by adjoint matrices in the basis $\{\de_\mu\}$
(equivalently, by transposed matrices, because all the matrices in questions
have real entries). Note that $H^*=H$.

Now identify $\Fun_0(\Y)$ and $\Sym$  (as vector spaces) via the correspondence
$\de_\mu\leftrightarrow s_\mu$. Then we may interpret $H^*=H$, $E^*$, and $F^*$
as operators in $\Sym$. {}From the definition of $F$ and the well--known
identity
$$
p_1 s_\la=\sum_{\la^\bullet\searrow\la}s_{\la^\bullet}
$$
it follows that the operator of multiplication by $p_1$ equals $-F^*$. Next,
the operator
$$
s_\mu\,\mapsto\, \sum_{\mu_\bullet\nearrow\mu}
(z)_{\mu/\mu_\bullet}(z')_{\mu/\mu_\bullet}s_{\mu_\bullet}
$$
equals $E^*$. Finally, recall the operator $G$ introduced just before Lemma 5.4
and note that
$$
G: s_\mu\mapsto|\mu|s_\mu.
$$
Therefore, one can write
$$
B_{z,z'}=-G(G-1+\th)-F^*E^*, \qquad \th:=zz'.
$$

We have to check that this operator commutes with $-F^*$, which is equivalent
to
$$
[G(G-1+\th)+EF, \, F]=0.
$$

Observe that $H=2G+\th$. Therefore
$$
G(G-1+\th)+EF=\frac14 H^2-\frac12 H+EF -\frac14\th^2+\frac12 \th.
$$
It is readily checked that this expression gives a central element in the
universal enveloping algebra of $\frak{sl}(2,\C)$. \qed
\enddemo

\head 6. The pre--generator $A_{z,z'}$ as a differential operator \endhead

Recall (see \S\S3.4--3.5) that the algebra $\La^\circ$ can be identified with
the polynomial ring $\R[q_1,q_2,\dots]$. Note that {\it any\/} linear operator
in the vector space of polynomials with countably many indeterminates can be
written as a differential operator, that is, as an infinite sum of differential
monomials with polynomial coefficients. Our aim here is to write in this form
the operator $A_{z,z'}:\bSym\to\bSym$ defined by \tht{5.1}.

\proclaim{Theorem 6.1} In the moment coordinates $q_1=p^\circ_2$,
$q_2=p^\circ_3,\dots$, the operator $A_{z,z'}$ defined by \tht{5.1} can be
written as the differential operator
$$
\multline A_{z,z'}=\sum_{i,j=1}^\infty
(i+1)(j+1)(q_{i+j}-q_iq_j)\frac{\pd^2}{\pd q_i\pd q_j}\\
-zz'\sum_{i=1}^\infty (i+1)q_i\frac{\pd}{\pd q_i}+ (z+z')\sum_{i=1}^\infty
(i+1)q_{i-1}\frac{\pd}{\pd q_i}\\+\sum_{i,j=0}^\infty(i+j+3)q_i
q_j\frac{\pd}{\pd q_{i+j+2}}-\sum_{i=1}^\infty (i+1)i q_i\frac{\pd}{\pd q_i}\,,
\endmultline \tag{6.1}
$$
where, by agreement, $q_0=1$.
\endproclaim

\demo{Proof} We will show that the operator $B_{z,z'}:\Sym\to\Sym$ (see
\tht{5.2}) can be written as the following differential operator in the
indeterminates $p_1,p_2,\dots$:
$$
\multline B_{z,z'}=\sum_{i,j=2}^\infty
ij(p_1p_{i+j-1}-p_ip_j)\frac{\pd^2}{\pd p_i\pd p_j}\\
-zz'\sum_{i=2}^\infty ip_i\frac{\pd}{\pd p_i}+ (z+z')\sum_{i=2}^\infty
ip_1p_{i-1}\frac{\pd}{\pd
p_i}\\+\sum_{i,j=1}^\infty(i+j+1)p_1p_ip_j\frac{\pd}{\pd p_{i+j+1}} -
\sum_{i=2}^\infty i(i-1)p_i\frac{\pd}{\pd p_i}
\endmultline \tag6.2
$$
Note that \tht{6.2} does not involve $\pd/\pd p_1$, so that we can reduce
\tht{6.2} modulo the relation $p_1=1$, by sending $p_i$ to $p^\circ_i=q_{i-1}$,
with the understanding that $p^\circ_1=q_0=1$. By virtue of Proposition 5.5
this will give us \tht{6.1}.

The main difficulty is to handle the sum in the right--hand side of \tht{5.2},
that is, the operator in $\Sym$ defined by
$$
s_\mu\mapsto\sum_{\mu_\bullet\nearrow\mu}
(z)_{\mu/\mu_\bullet}(z')_{\mu/\mu_\bullet}s_{\mu_\bullet}\,.
$$
Given a box $\square$ with the row coordinate $i$ and the column coordinate
$j$, we denote by $c(\square)$ the content $j-i$. The latter operator can be
written as a linear combination of three operators,
$$
C'_2+(z+z')C'_1+zz'C'_0,
$$
where by $C'_k:\Sym\to\Sym$ we denote the operator
$$
C'_k: s_\mu\mapsto
\sum_{\mu_\bullet\nearrow\mu}(c(\mu/\mu_\bullet))^ks_{\mu_\bullet}\,, \qquad
k=0,1,2,\dots
$$
In this notation the operator $B_{z,z'}$ is written as follows
$$
B_{z,z'}=-G(G-\bold1+zz'\bold1)+p_1(C'_2+(z+z')C'_1+zz'C'_0). \tag6.3
$$

By the very definition of $G$ we have
$$
G=\sum_{i=1}^\infty ip_i\frac{\pd}{\pd p_i}\,.\tag6.4
$$
The operator $C'_0$ is also easy to write:
$$
C'_0=\pd/\pd p_1.\tag6.5
$$
Indeed, to see this, one can use the fact that $\pd/\pd p_1$ is adjoint to the
operator of multiplication by $p_1$, which has the form
$$
s_\mu\mapsto \sum_{\mu^\bullet\searrow\mu}s_{\mu^\bullet}\,.
$$
The operator $C'_1$ and especially the operator $C'_2$ are more involved.

To handle them it is convenient to introduce auxiliary  operators
$C_k:\Sym\to\Sym$ by
$$
C_k: s_\mu\mapsto \left(\sum_{\square\in\mu}(c(\square))^k\right)s_\mu\,,
$$
summed over all boxes contained in $\mu$, and observe that $C'_0=\pd/\pd p_1$
implies
$$
C'_k=\left[\frac{\pd}{\pd p_1}\,,\; C_k\right].
$$

We will employ the following results proved by Lascoux and Thibon in \cite{LT,
Proposition 3.3}:

\proclaim{Lemma 6.2} Let $t$ and $u$ be formal variables. The exponential
generating series for the operators $C_k$ has the form
$$
\sum_{k=1}^\infty C_k \frac{t^k}{k!}=\frac{V_0-1}{(e^t-1)(1-e^{-t})} - G,
$$
where $V_0$ is the constant term of the ``vertex operator''
$$
\multline V(u)=\sum_{m=-\infty}^{+\infty} V_mu^{-m}\\
:=\exp\left(\sum_{k=1}^\infty(e^{kt}-1)\frac{u^k}{k}p_k\right)
\exp\left(\sum_{k=1}^\infty(1-e^{-kt})u^{-k}\frac{\pd}{\pd p_k}\right).
\endmultline
$$

Moreover, the following commutation relations hold:
$$
\left[V_m\,,\; k\frac{\pd}{\pd p_k}\right]=(1-e^{kt})V_{m+k}.
$$
In particular,
$$
\left[\frac{\pd}{\pd p_1}\,,\;V_0\right]=(e^t-1)V_1.
$$
\endproclaim

{}From Lemma 6.2 we deduce explicit expressions for the operators $C'_1$ and
$C'_2$:

\proclaim{Lemma 6.3} We have
$$
\gather
 C'_1=\sum_{i=1}^\infty (i+1)p_i\frac{\pd}{\pd p_{i+1}} \tag6.6\\
C'_2=\sum_{i,j=1}^\infty ijp_{i+j-1}\frac{\pd}{\pd p_i} \frac{\pd}{\pd
p_j}+\sum_{i,j=1}^\infty(i+j+1)p_ip_j\frac{\pd}{\pd p_{i+j+1}}. \tag6.7
\endgather
$$
\endproclaim

Note that this particular  result can also be obtained by the more elementary
approach used in \cite{FW}.

\demo{Proof of Lemma 6.3} {}From Lemma 6.2 we get
$$
\left[\frac{\pd}{\pd p_1}\,,\;C_k\right]=\text{\rm coefficient of $t^ku^{-1}$
in\; $k!\,\frac{V(u)}{1-e^{-t}}$}\,, \quad k=1,2,\dots
$$
Next, write
$$
V(u)=\exp\left(\sum_{r=1}^\infty a_rt^r\right)\exp\left(\sum_{r=1}^\infty
b_rt^r\right),
$$
where
$$
a_r=\sum_{k=1}^\infty\frac{k^ru^k}{r!k}p_k\,, \qquad
b_r=(-1)^{r-1}\sum_{k=1}^\infty\frac{k^ru^{-k}}{r!}\frac{\pd}{\pd p_k}\,.
\tag6.8
$$
Then
$$
\multline C'_1=\left[\frac{\pd}{\pd p_1}\,,\;C_1\right]=\text{the coefficient
of
$t^2u^{-1}$\; in}\\
\frac{\exp(a_1t+a_2t^2+\dots)\exp(b_1t+b_2t^2+\dots)}{1-\frac12 t+\frac16
t^2+\dots}
\endmultline \tag6.9
$$
and
$$
\multline C'_2=\left[\frac{\pd}{\pd p_1}\,,\;C_2\right]=\text{the coefficient
of
$t^3u^{-1}$\; in}\\
2\cdot\frac{\exp(a_1t+a_2t^2+a_3t^3+\dots)\exp(b_1t+b_2t^2+b_3t^3+\dots)}{1-\frac12
t+\frac16 t^2-\frac1{24}t^3+\dots}\,,
\endmultline \tag6.10
$$
where the dependence in $u$ is hidden in \tht{6.8}. First, we compute the
coefficients of $t^2$ and $t^3$ in \tht{6.9} and \tht{6.10}, respectively.

The coefficient of $t^2$ in \tht{6.9} equals
$$
a_2+\frac12 a_1^2+a_1b_1+\frac12a_1+b_2+\frac12b_1^2+\frac12b_1+\frac1{12}.
$$
{}From \tht{6.8} it immediately follows that a nonzero contribution
to the term with $u^{-1}$ can come from $a_1b_1$, $b_2$, and
$\frac12b_1$ only. It turns out that the total contribution of
$b_2+\frac12b_1$ equals 0. This gives \tht{6.6}.

Likewise, the coefficient of $t^3$ in \tht{6.10} equals
$$
\gather 2a_3+a_2+2a_1a_2+2a_2b_1+\frac12a_1^2+\frac13a_1^3+a_1^2b_1+2a_1b_2
+a_1b_1\\
+a_1b_1^2+\frac16a_1+2b_3+b_2+2b_1b_2+\frac12b_1^2+\frac13b_1^3+\frac16b_1
\endgather
$$
and a nonzero contribution to the term with $u^{-1}$ can come from
$$
2a_2b_1+a_1^2b_1+2a_1b_2+a_1b_1+a_1b_1^2+2b_3+b_2+\frac16b_1
$$
only. It turns out that the total contribution of $a_1b_1+2a_2b_1+2a_1b_2$
equals 0, and the same holds for $2b_3+b_2+\frac16b_1$. Finally, we get
\tht{6.7}. \qed
\enddemo

Substituting the expressions \tht{6.4}, \tht{6.5}, \tht{6.6}, and \tht{6.7}
into \tht{6.3} we get after cancellations the desired expression \tht{6.2}.
\qed
\enddemo

Let $\bar A_{z,z'}$ denote the closure of the operator $A_{z,z'}$ in $C(\Om)$
(recall that the closure exists according to Proposition 1.4 (1)). The next
result will be used in the proof of Theorem 7.1 below.

\proclaim{Corollary 6.4} Let $f(\om)$ be a smooth cylinder function in the
moment coordinates, that is, $f(\om)=g(q_1(\om),\dots,q_m(\om))$ for a certain
$m=1,2,\dots$ and a certain smooth function $g(q_1,\dots,q_m)$ in $m$ real
variables, in a neighborhood of $[-1,1]^m$. Then $f$ enters the domain of $\bar
A_{z,z'}$. Moreover, $\bar A_{z,z'} f$ is also a cylinder function, which can
be obtained via application of the suitably truncated differential expression
\tht{6.1} to the function $g$.
\endproclaim

\demo{Proof} First of all, note that here ``truncation'' means that we keep in
\tht{6.1} only terms not containing derivatives $\pd/\pd q_i$ with $i>m$. It is
worth noting that the resulting cylinder function depends on the larger number
of variables, $2m$, because of the presence of the variables $q_{k+l}$ in
\tht{6.1}. However, this does not cause problems.

To prove the claim of the proposition we observe that the function
$g(q_1,\dots,q_m)$, together with its partial derivatives of up to second
order, can be approximated by a sequence $\{g_n(q_1,\dots,q_m)\}$ of
polynomials, uniformly on the cube $[-1,1]^m\subset\R^m$. Let $[A_{z,z'}]$
stand for the truncated differential operator as explained above. The
application of $\bar A_{z,z'}$ to the function $\om\mapsto
g_n(q_1(\om),\dots,q_m(\om))$ is reduced to the application of $[A_{z,z'}]$ to
$g_n$. Since $g_n\to g$ and $[A_{z,z'}]g_n\to[A_{z,z'}]g$ uniformly on the cube
$[-1,1]^{2m}$, we see that $f$ belongs to the domain of $\bar A_{z,z'}$ and
$\bar A_{z,z'} f$ is given by the cylinder function $[A_{z,z'}]g$. \qed
\enddemo

\head 7. The limit process\endhead

The theorems of this section are almost direct consequences of the results
established in sections 1, 5, and 6. Note that application of Propositions 1.6
and 1.7 is justified, because their hypothesis is satisfied due to the theorem
of \S3.6.

\proclaim{Theorem 7.1} The Markov semigroup $\{T(t)\}$ in $C(\Om)$ constructed
in \S5 gives rise to a diffusion process $\ome_{z,z'}(t)$ in $\Om$.
\endproclaim

By a diffusion process we mean a strong Markov process (which can start from
any point or any probability distribution) with continuous sample paths.

\demo{Proof} Once the existence of a Markov semigroup $\{T(t)\}$ is established
(Theorem 5.1), the existence of the  corresponding strong Markov process is
guaranteed by a general theorem stated above as Proposition 1.5. It remains to
prove that the sample paths are continuous almost surely.

Since our semigroup acts in the space of continuous functions on a compact
space, the application of Riesz's theorem implies the existence of a transition
function (see, e.g., \cite{L, \S7.7, Thm. 1}). The continuity property holds if
the transition function obeys the Dynkin--Kinney condition (\cite{W, \S10.3} or
\cite{EK2, Ch. 4, \tht{2.35}}). This condition in turn holds if for any point
$\om\in\Om$ and any its neighborhood $U\ni\om$ one can find a function $f\in
C(\Om)$ with the following properties (see \cite{EK2, Ch. 4, Remark 2.10}):

(1) $f$ is contained in $D(\bar A_{z,z'})$, the domain of the generator $\bar
A_{z,z'}$ of the semigroup;

(2) $\bar A_{z,z'}f(\om)=0$;

(3) $\Vert f\Vert=f(\om)$ and the supremum of $f$ outside $U$ is strictly less
than $f(\om)$.

Such functions can be built using Corollary 6.4. Indeed, take $m$ so large that
one can find inside $U$ a neighborhood of the form
$$
\{\om'\in\Om:\quad |q_i(\om')-q_i(\om)|<\epsi, \quad 1\le i\le m\}.
$$
Then take as $f$ a cylinder function as in Corollary 6.4, where $g$ equals 1 in
a very small neighborhood of the point $(q_1(\om),\dots,q_m(\om))\in[-1,1]^m$
and then rapidly abates to 0. Since the differential operator \tht{6.1} does
not have a constant term, $\bar A_{z,z'}f$ vanishes about $\om$. \qed
\enddemo

Let $P_{z,z'}$ be the boundary measure on $\Om$ corresponding to the coherent
system with parameters $(z,z')$, and consider the inner product in
$\Sym^\circ\subset C(\Om)$ given by
$$
(f,g)_{z,z'}=\langle f\cdot g\rangle_{P_{z,z'}}:=\int_\Om
f(\om)g(\om)P_{z,z'}(d\om).
$$

\proclaim{Theorem 7.2} \tht{1} The space $\Sym^\circ$ can be decomposed into a
direct sum of eigenspaces of the pre--generator $A_{z,z'}$, and this
decomposition is orthogonal with respect to the above inner product.

\tht{2} The spectrum of $A_{z,z'}$ is $\{0\}\cup\{-\si_m: m=2,3,\dots\}$ where
$$
\si_m=m(m-1+zz'), \qquad m=2,3,\dots
$$

\tht{3} The eigenvalue\/ $0$ is simple, and the multiplicity of $-\si_m$ equals
the number of partitions of $m$ without parts equal to $1$.
\endproclaim

\demo{Proof} (1) This is a fact of linear algebra because the pre--generator
$A_{z,z'}$ is symmetric (Proposition 1.7) and preserves the filtration of
$\Sym^\circ$.

(2) Let $I$ denote the principal ideal in $\Sym$ generated by $p_1-1$. We have
$$
\Sym=\R[p_1,p_2,p_3,\dots]=\R[p_2,p_3,\dots]\oplus I
$$
so that we may identify $\Sym^\circ$ with $\R[p_2,p_3,\dots]$. It follows from
\tht{5.2} that for any homogeneous element $f\in\Sym$ of degree $m$,
$$
B_{z,z'}f=-m(m-1+zz')f+g+h, \qquad \text{\rm where $g\in I$ and $\deg h<m$}.
$$
In particular, this is true for any monomial $p_2^{m_2}p_3^{m_3}\dots$ from
$\R[p_2,p_3,\dots]$. Taking into account Proposition 5.5 we conclude that the
spectrum of $A_{z,z'}$ is as indicated in claim (2), and the multiplicity of
$-\si_m$ equals the number of solutions in nonnegative integers of the equation
$$
2m_2+3m_3+\dots=m,
$$
which proves claim (3). \qed
\enddemo

\proclaim{Theorem 7.3} \tht{1} The process $\ome_{z,z'}(t)$ constructed in
Theorem\/ {\rm 7.1} has the boundary measure $P_{z,z'}$ as a unique stationary
distribution.

\tht{2} It is also a symmetrizing measure.

\tht{3} The process is ergodic in the sense that for any $f\in C(\Om)$,
$$
\lim_{t\to+\infty}\Vert T(t)f-\langle f\rangle_{P_{z,z'}}\cdot1\Vert=0,
$$
where $\Vert\,\cdot\,\Vert$ is the norm of the Banach space $C(\Om)$ and\/ $1$
is the constant function equal to one.
\endproclaim

\demo{Proof} Consider the orthogonal decomposition of $\Sym^\circ$ onto
eigenspaces afforded by Theorem 7.2:
$$
\Sym^\circ=\R1\oplus\bigoplus_{m=2}^\infty\Sym^\circ_m\,. \tag7.1
$$
The operator $T(t)$ leaves invariant the constant 1 and acts in $\Sym^\circ_m$
as multiplication by $\exp(-\si_mt)$. Note that the direct sum decomposition in
\tht{7.1} is understood in purely algebraic sense: for any vector
$f\in\Sym^\circ$, its expansion $f=f_0+f_2+f_3+\dots$ has finitely many nonzero
components.

(1) The fact that $P_{z,z'}$ is an invariant distribution follows from
Proposition 1.6. To prove uniqueness we observe that if $P$ is an invariant
distribution then $\langle f\rangle_P=0$ for any $f\in\Sym^\circ_m$,
$m=2,3,\dots$. Therefore, for any $f=f_0+f_2+f_3+\dots\in\Sym^\circ$ with
$f_0=c1$ we have
$$
\langle f\rangle_P=\langle f_0\rangle_P=c. \tag7.2
$$
Since $\Sym^\circ$ is dense in $C(\Om)$, $P$ is determined uniquely.

(2) The claim to be proved is equivalent to the fact that
$$
\int_\Om(\bar A_{z,z'}f(\om))g(\om)P_{z,z'}(d\om)
$$
is symmetric under transposition $f\leftrightarrow g$, for any $f$ and $g$ in
the domain of $\bar A_{z,z'}$. It suffices to check this for $f$ and $g$ in
$\Sym^\circ\subset C(\Om)$ and with $A_{z,z'}$ replacing $\bar A_{z,z'}$, which
follows from Proposition 1.7.

(3) {}From the decomposition \tht{7.1} it is evident that for any
$f\in\Sym^\circ$, we have $T(t)f\to f_0$ as $t\to+\infty$. Here the convergence
holds in a finite--dimensional space invariant under the semigroup $T(t)$,
hence $T(t)f$ converges to $f_0$ in norm, too. Together with \tht{7.2} this
proves claim (3) for $f\in\Sym^\circ$. Then it is evident that it also holds
for any $f\in C(\Om)$. \qed
\enddemo

Recall that any function $f(\om)$ from $\Sym^\circ\subset C(\Om)$ can also be
viewed  as a polynomial in the moment coordinates $q_1=q_1(\om)$,
$q_2=q_2(\om)$, \dots. With this understanding, we set, for any two functions
$f,g\in\Sym^\circ$
$$
\Ga(f,g)(\om)=\sum_{i,j=1}^\infty \Ga_{ij}(\om) \frac{\pd f}{\pd
q_i}(\om)\frac{\pd g}{\pd q_j}(\om)\,, \tag7.3
$$
where
$$
\Ga_{ij}(\om)=(i+1)(j+1)(q_{i+j}(\om)-q_i(\om)q_j(\om)), \qquad i,j=1,2,\dots
\tag7.4
$$
Note that the sum is actually finite, because the partial derivatives with
sufficiently large indices vanish.

\proclaim{Theorem 7.4} For any $f,g\in\Sym^\circ\subset C(\Om)$
$$
-\int_\Om
A_{z,z'}f(\om)g(\om)P_{z,z'}(d\om)=\int_{\Om}\Ga(f,g)(\om)P_{z,z'}(d\om)
$$
\endproclaim

The point here is that both the pre--generator $A_{z,z'}$ and the boundary
measure $P_{z,z'}$ depend on the parameters $(z,z')$ while $\Ga(f,g)$ does not:
in the right--hand side, the parameters enter $P_{z,z'}$ only.

\demo{Proof} We abbreviate $A=A_{z,z'}$ and
$\langle\,\cdot\,\rangle=\langle\,\cdot\,\rangle_{P_{z,z'}}$. Let us show that
$$
2\Ga(f,g)=A(fg)-(Af)g-f(Ag), \qquad f,g\in\Sym^\circ. \tag7.5
$$
Indeed, by Theorem 6.1, $A$ is a second order differential operator in the
moment coordinates. Therefore, its first order terms do not contribute to the
right--hand side of \tht{7.5}. Writing
$$
A=\sum_{i,j=1}^\infty \Ga_{ij}\frac{\pd^2}{\pd q_i\pd q_j}\, +\,\text{first
order terms}
$$
we see that the right--hand side of \tht{7.5} is equal to
$$
2\sum_{i,j=1}^\infty \Ga_{ij}\frac{\pd f}{\pd q_i}\frac{\pd g}{\pd q_j}\,,
$$
which is precisely the definition of $\Ga(f,g)$.

Next, as we already pointed out above, the expectation
$\langle\,\cdot\,\rangle$ vanishes on all eigenspaces of $A$, except that
corresponding to the eigenvalue 0. Consequently $\langle\,\cdot\,\rangle$
vanishes on the range of the operator $A$. Applying the expectation to the both
sides of \tht{7.5} and using the fact that $A$ is symmetric, we get the desired
formula. \qed
\enddemo

\example{Remark 7.5} Note that for any $\om\in\Om$, the infinite matrix
$[\Ga_{ij}(\om)]_{i,j=1}^\infty$ is nonnegative definite. Indeed, recall that
$q_i=q_i(\om)$ is the $i$th moment of the Thoma measure $\nu_\om$ on $[-1,1]$.
It follows that for any sequence $c_1,c_2,\dots$ of real numbers with finitely
many nonzero entries, the quadratic form $\sum_{i,j}(q_{i+j}-q_iq_j)c_ic_j$
equals the variance of the function $\sum_i c_ix^i$ with respect to $\nu_\om$.
Cf. Schmuland \cite{S, p.~255}. It is tempting to regard the quantity
$\Ga(f,f)$ defined by \tht{7.3}--\tht{7.4} as a square field (carr\'e du
champs). Similar expressions already appeared in works on measure--valued
diffusions, see, e.g., Overbeck--R\"ockner--Schmuland \cite{OvRS}.
\endexample

\bigskip

{\smc A.~Borodin}: Mathematics 253-37, Caltech, Pasadena, CA 91125, U.S.A.

\medskip

E-mail address: {\tt borodin\@caltech.edu}

\bigskip

{\smc G.~Olshanski}: Dobrushin Mathematics Laboratory, Institute for
Information Transmission Problems, Bolshoy Karetny 19, 127994 Moscow GSP-4,
RUSSIA.

\medskip

E-mail address: {\tt olsh\@online.ru}


\Refs

\widestnumber\key{GNW2}

\ref\key B \by A. Borodin \paper Harmonic analysis on the infinite symmetric
group and the Whittaker kernel \jour St.~Petersburg Math. J. \vol 12 \yr 2001
\issue 5 \pages 733-759
\endref

\ref\key BO1 \by A.~Borodin and G.~Olshanski \paper Point processes and the
infinite symmetric group \jour Math. Research Lett. \vol 5 \yr 1998 \pages
799--816; arXiv: math.RT/9810015
\endref

\ref\key BO2 \by A.~Borodin and G.~Olshanski \paper Distributions on
partitions, point processes and the hypergeometric kernel \jour Comm. Math.
Phys. \vol 211 \yr 2000 \pages 335--358;  arXiv: math.RT/9904010
\endref

\ref\key BO3 \by A.~Borodin and G.~Olshanski \paper Harmonic functions on
multiplicative graphs and interpolation polynomials \jour Electronic J. Comb.
\vol 7 \yr 2000 \pages paper \#R28; arXiv: math/9912124
\endref

\ref\key BO4 \by A.~Borodin and G.~Olshanski \paper Z--Measures on partitions,
Robinson--Schensted--Knuth correspondence, and $\beta=2$ random matrix
ensembles \inbook in: Random matrix models and their applications (P.~M.~Bleher
and A.~R.~Its, eds). Mathematical Sciences Research Institute Publications {\bf
40} \publ Cambridge Univ. Press \yr 2001 \pages 71--94; arXiv: math/9905189
\endref

\ref\key BO5 \by A.~Borodin and G.~Olshanski \paper Random partitions and the
Gamma kernel \jour Advances in Math. \vol 194 \yr 2005 \issue 1 \pages
141--202; arXiv: math-ph/0305043
\endref

\ref\key BO6 \by A.~Borodin and G.~Olshanski \paper Markov processes on
partitions \jour Prob. Theory and Related Fields \vol 135 \yr 2006 \issue 1
\pages 84--152; arXiv: math-ph/0409075
\endref

\ref\key BO7 \by A.~Borodin and G.~Olshanski \paper Meixner polynomials and
random partitions \jour Moscow Math. J. \vol 6 \yr 2006 \issue 4 \pages
629--655; arXiv: math.PR/0609806
\endref

\ref\key Ed \by A. Edrei \paper On the generating functions of totally positive
sequences {\rm II} \jour J. Analyse Math. \vol 2 \yr 1952 \pages 104--109
\endref

\ref\key EK1 \by S.~N.~Ethier and T.~G.~Kurtz \paper The
infinitely--many--neutral--alleles diffusion model \jour Adv. Appl. Prob. \vol
13\yr 1981 \pages 429--452 \endref

\ref\key EK2 \by S.~N.~Ethier and T.~G.~Kurtz \book Markov processes --
Characterization and convergence \publ Wiley--Interscience \publaddr New York
\yr 1986 \endref

\ref\key FW \by I.~B.~Frenkel and W.~Wang \paper Virasoro algebra and wreath
product convolution \jour J. Algebra \vol 242 \yr 2001\issue 2 \pages 656--671
\endref

\ref\key F \by J.~Fulman \paper Stein's method and Plancherel measure of the
symmetric group \jour Trans. Amer. Math. Soc. 357 \yr 2005 \issue 2\pages
555-570; arXiv: math.RT/0305423
\endref

\ref\key IO \by V.~Ivanov and G.~Olshanski \paper Kerov's central limit theorem
for the Plancherel measure on Young diagrams \inbook In: Symmetric functions
2001. Surveys of developments and perspectives. Proc. NATO Advanced Study
Institute (S.~Fomin, editor), Kluwer, 2002, pp. 93--151; arXiv: math/0304010
\endref

\ref\key KSK \by J.~G.~Kemeny, J.~L.~Snell, and A.~W.~Knapp \book Denumerable
Markov chains \publ Springer \publaddr NY \yr 1976
\endref

\ref\key K \by S.~V.~Kerov \book Asymptotic representation theory of the
symmetric group and its applications in analysis \publ Amer. Math. Soc.,
Providence, RI, 2003, 201 pp
\endref

\ref\key KOO \by S.~Kerov, A.~Okounkov, G.~Olshanski \paper The boundary of
Young graph with Jack edge multiplicities \jour Intern. Math. Res. Notices \yr
1998 \issue 4 \pages 173--199; arXiv: q-alg/9703037
\endref

\ref \key KO \by S.~Kerov and G.~Olshanski \paper Polynomial functions on the
set of Young diagrams \jour Comptes Rendus Acad.\ Sci.\ Paris S\'er. I \vol 319
\yr 1994 \pages 121--126
\endref

\ref \key KOV1 \by S.~Kerov, G.~Olshanski, and A.~Vershik \paper Harmonic
analysis on the infinite symmetric group. A deformation of the regular
representation \jour Comptes Rendus Acad. Sci. Paris, S\'er. I \vol 316 \yr
1993 \pages 773--778
\endref

\ref\key KOV2 \by S.~Kerov, G.~Olshanski, and A.~Vershik \paper Harmonic
analysis on the infinite symmetric group \jour Invent. Math. \vol 158 \yr 2004
\pages 551--642;  arXiv: math.RT/0312270
\endref

\ref\key L \by J.~Lamperti \book Stochastic processes, A survey of the
mathematical theory \publ Springer \yr 1977
\endref

\ref\key LT \by A.~Lascoux and J.--Y.~Thibon \paper Vertex operators and the
class algebras of symmetric groups \jour  J. Math. Sci. (N. Y.)  \vol 121 \yr
2004\issue 3 \pages 2380--2392; arXiv: math/0102041
\endref

\ref\key Ma \by I.~G.~Macdonald \book Symmetric functions and Hall polynomials
\bookinfo 2nd edition \publ Oxford University Press \yr 1995
\endref

\ref\key Ok \by A.~Okounkov \paper $SL(2)$ and $z$--measures \inbook in: Random
matrix models and their applications (P.~M.~Bleher and A.~R.~Its, eds).
Mathematical Sciences Research Institute Publications {\bf 40} \publ Cambridge
Univ. Press \yr 2001 \pages 407--420; arXiv: math.RT/0002136
\endref

\ref \key OO \by A.~Okounkov and G.~Olshanski \paper Shifted Schur functions
\jour Algebra i Analiz \vol 9 \issue 2 \yr 1997 \pages 73--146 \lang Russian
\transl English translation: St.~Petersburg Math. J. {\bf 9} (1998), no.~2,
239--300; arXiv: q-alg/9605042
\endref

\ref\key Ol1 \by G. Olshanski \paper Point processes related to the infinite
symmetric group \inbook In: The orbit method in geometry and physics: in honor
of A.~A.~Kirillov (Ch.~Duval et al., eds.), Progress in Mathematics {\bf 213},
Birkh\"auser, 2003, pp. 349--393;  arXiv: math.RT/9804086
\endref

\ref\key Ol2 \by G. Olshanski \paper An introduction to harmonic analysis on
the infinite symmetric group \inbook In: Asymptotic combinatorics with
applications to mathematical physics \ed A.~M.~ Vershik \bookinfo A European
mathematical summer school held at the Euler Institute, St.~Petersburg, Russia,
July 9--20, 2001 \publ Springer Lect. Notes Math. {\bf 1815}, 2003, 127--160;
arXiv: math.RT/0311369
\endref

\ref \key OlRV \by G.~Olshanski, A.~Regev, and A.~Vershik \paper
Frobenius--Schur functions \inbook Studies in memory of Issai Schur (A.~Joseph,
A.~Melnikov, R.~Rentschler, eds), Progress in Mathematics {\bf 210},
Birkh\"auser, 2003, pp. 251--300; arXiv: math/0110077
\endref

\ref \key OvRS \by L.~Overbeck, M.~R\"ockner, and B.~Schmuland \paper An
analytic approach to Fleming--Viot processes with interactive selection \jour
Ann. ~Prob. \vol 23 \yr 1995 \pages 1--36
\endref

\ref\key P \by L.~A.~Petrov \paper Two--parameter family of diffusion processes
in the Kingman simplex \paperinfo Preprint 2007,  arXiv:0708.1930 [math.PR]
\endref

\ref\key S \by B.~Schmuland \paper A result on the infinitely many neutral
alleles diffusion model\jour J. Appl. Prob. \vol 28\yr 1991 \pages 253--267
\endref

\ref \key T \by E.~Thoma \paper Die unzerlegbaren, positive-definiten
Klassenfunktionen der abz\"ahlbar unendlichen, symmetrischen Gruppe \jour
Math.~Zeitschr. \vol 85 \yr 1964 \pages 40--61
\endref

\ref\key VK \by A.~M.~Vershik and S.~V.~Kerov \paper Asymptotic theory of
characters of the symmetric group \jour Funct. Anal. Appl. \vol 15 \yr 1981
\pages 246--255
\endref

\ref\key W \by A.~D.~Wentzell \book A course in the theory of stochastic
processes \publ McGraw-Hill International \publaddr New York \yr 1981 \endref

\endRefs
\enddocument
\end